\RequirePackage{fix-cm}
\documentclass[smallextended]{svjour3}       % onecolumn (second format)
\smartqed  % flush right qed marks, e.g. at end of proof

\usepackage{graphicx}
\usepackage{amsmath}         % AMSLaTeXºê°ü ÓÃÀ?ÅÅ³ö?ü?ÓÆ¯ÁÁµÄ¹«Ê?
\usepackage{amssymb}
\usepackage{mathrsfs}
\usepackage{appendix}
\usepackage[latin1]{inputenc}
\newcommand \mH{\mathcal{H}}
\newcommand \mM{\mathcal{M}}
\newcommand \mG{\mathcal{G}}
\newcommand \mF {\mathcal{F}}

\newcommand \mO {\mathcal{O}}

\newcommand \mA {\mathcal{A}}
\newcommand \mC {\mathcal{C}}

\newcommand \pr {\mathbb{P}}
\newcommand \E {\mathbb{E}}

\newcommand \R {\mathbb{R}}

\newcommand \la {\lambda}
\newcommand \alp {\alpha}

\newcommand \kap {\kappa}
\newcommand \sig {\sigma}

\newtheorem{assumption}{Assumption}%[section]
\newcommand \hU {\hat{U}}

\begin{document}

\title{Optimal  per-loss reinsurance  and investment to minimize the probability of drawdown}

\author{Xia Han  \and Zhibin Liang }


\institute{Xia Han \at
             School of Mathematical Sciences, Nanjing Normal University,\\
            Jiangsu 210023, P.R.China \\
             307717378@qq.com
           \and
          Zhibin Liang,  Corresponding author  \at
             School of Mathematical Sciences, Nanjing Normal University,\\
             Jiangsu 210023, P.R.China \\
           liangzhibin111@hotmail.com}

\date{Received: date / Accepted: date}

\maketitle

\begin{abstract}
In this paper, we study an optimal reinsurance-investment  problem in a risk model with two dependent classes of insurance business, where the two claim number processes are correlated through a common shock component. We assume that the insurer can purchase per-loss reinsurance for each line of business and invest its surplus in a financial market consisting of a risk-free asset  and a risky asset. Under the criterion of minimizing the probability of drawdown, the closed-form expressions of the optimal reinsurance-investment strategy and the corresponding  value function are obtained. We show that the optimal reinsurance strategy is in the form of pure excess-of-loss reinsurance strategy under the expected value  principle, and under the variance premium principle, the optimal reinsurance strategy is in the form of  pure quota-share reinsurance. Furthermore, we extend our model to the   case where the  insurance company involves $n$ $(n\geq3)$ dependent classes of insurance business and the optimal results are derived explicitly as well.

%\subclass{   62P05 \and 90C39 \and 91B30 \and 93E20  }
\end{abstract}

\section{Introduction}\label{sec:1}
In the past decades, risk models taking reinsurance  into consideration have received a great deal of attention in the literature. With reinsurance, insurers  can transfer some of their risks to other insurers    at the expense of making less potential profit, and hence finding optimal reinsurance strategy to balance their risks and profit is of great interest to them.  Subject to controlling reinsurance with or without controlling investment, optimization problems such as minimizing the probability of ruin,  maximizing the expected utility of terminal wealth,  as well as the traditional mean-variance  optimization criterion  have become a popular research topic in the actuarial literature.  The technique of stochastic control theory and the corresponding Hamilton- Jacobi-Bellman (HJB) equation are widely used to cope with these problems.

Much of  the existing literature  considered  problems under which  the insurance company is constrained to buy either pure quota-share reinsurance, excess-of-loss reinsurance, or a combination of the two. Promislow and Young \cite{PY05}  found the optimal investment and quota-share reinsurance strategies to minimize the probability of ruin of an insurer who faces a claim process that  follows Brownian motion with drift.  Liang and Guo \cite{LG11} maximized expected exponential utility of terminal wealth by finding the optimal combination of quota-share and excess-of-loss reinsurance. Under the mean-variance criterion,  Li et al. \cite {LRZ17} derived an equilibrium reinsurance-investment strategy in a financial market where the  surplus process is assumed to follow the classical Cram\'er-Lundberg model. Some other researchers have found optimal reinsurance strategies for various optimization problems without restricting the form of the reinsurance. For example,  %Bai et al. \cite{BCZ13} considered the two-dimensional reinsurance policy in a dynamic setting and showed that a two-dimensional excess-of-loss reinsurance policy is an optimal form that minimizes the ruin probability of the controlled diffusion process.
 Zhang et al. \cite{ZMZ16} investigated optimal investment and reinsurance problems in which  the insurer purchases a general reinsurance policy and the reinsurance is priced according to the mean-variance premium principle.   Liang and Young \cite{LY18} computed the optimal investment and per-loss reinsurance strategies for an insurance company facing a compound Poisson claim process; they assumed that the reinsurer used an expected-value premium principle.  Han et al. \cite{HLY18}  found the general optimal reinsurance strategy  to minimize the probability of drawdown under the mean-variance premium principle and observed  that the optimal reinsurance strategy is identical to the one for minimizing the probability of ruin. 

 As we know,  drawdown occurs when the value of surplus drops to a fixed proportion of its maximum value. 
 Portfolio optimization problems related to drawdown risks have long  focused on maximizing  the long-term growth rate of a portfolio subject to a strict drawdown constraint;  see  Grossman and Zhou \cite{GZ93}, Cvitani\'c and Karatzas \cite{CK95} for  example. However, if the individual  is consuming continuously from its investment account %or the insurance company keeps paying expensive reinsurance premium,
   then one cannot prevent drawdown, so minimizing the probability of drawdown is a reasonable, objective goal. Angoshtari et al. \cite{ABY16b} and Chen et al. \cite{CLLL15} computed the optimal investment strategy to minimize the probability that drawdown occurs during the investor's life.  Angoshtari et al. \cite{ABY16a} and Han et al.  \cite {HLY19}\cite{HLY18} minimized the probability of drawdown over an infinite-time horizon and discovered  that the strategy which minimizes the probability of ruin also minimizes the probability of drawdown.  

In this paper, we  study the optimization  problem  for an insurer who wishes to minimize the probability of drawdown. We suppose that the insurer has  two dependent classes of insurance business, where the two claim number processes are correlated through a common shock component; and  the insurer can purchase per-loss  reinsurance and invest its surplus in a financial market consisting of one risky asset and one risk-free asset. To make the optimization  problem tractable and to obtain the explicit solutions,  we consider the diffusion approximation model of the  Cram\'er-Lundberg model.   Although the problem of optimal reinsurance has been widely studied, to the best of our knowledge,  only  Bai et al. \cite{BCZ13} investigated  the minimum ruin probability without restricting the form of reinsurance   under the   dependence structure.  According to the expected value  principle, they derived the explicit expression of   the optimal reinsurance strategy, which is independent of the surplus process; and showed that the excess-of-loss reinsurance is an optimal form. 
Compared to  Bai et al. \cite{BCZ13}, there are three main differences and contributions  in our paper. Firstly, we consider a more general criterion that minimizes the probability of drawdown, which includes minimizing the probability of ruin as a special case.  Since the domain and boundary conditions for drawdown problem  are very different  from  ruin, a verification theorem is  proved to show that a $\mC^{2,1}$ solution of the HJB equation coincides with the minimum drawdown function.     Secondly, while purchasing  reinsurance to reduce  risk exposure,  we also allow the insurer to  invest the  surplus into the risky and  risky-free assets to   increase its profit. The resulting  optimal expressions    strongly depend on the  values of the wealth and the drawdown level, which makes the drawdown problem more challenging and practical.  Thirdly, the optimization   problem  is considered not only for the expected value principle but also for the variance premium principle; and we find that the optimal strategies under the two  premium principles are totally different. Interestingly, when the same safety loading applies to all classes under the variance premium principle, we work out an optimal retention level which holds for all classes, and falls into the interval [0, 1].   

By the technique of stochastic dynamic programming, we obtain the closed-form expressions of optimal reinsurance-investment strategy and the corresponding  value function. We show that the optimal reinsurance strategy is in the form of pure excess-of-loss reinsurance strategy under the expected value principle, and under the variance premium principle, the optimal reinsurance strategy is in the form of  pure quota-share reinsurance.  These results are expected from  the work of Hipp and Taksar \cite{HT10}, in which   one  class of  insurance business was considered.   Furthermore, we extend our model to the   case where the  insurance company involves $n$ $(n\geq3)$ dependent classes of insurance business and the optimal results are derived  explicitly as well. 

The rest of the paper is organized as follows.  In Section \ref{sec:2}, we present the model and optimization problem.  In Section \ref{sec:3}, we prove a verification theorem, which is used to find our value function. The explicit expressions for the optimal strategy and the corresponding value function are derived in Section \ref{sec:4}.   In Section \ref{sec:5}, we further explore the optimization problem when  the insurance company involves $n$ $(n\geq3)$ dependent classes of insurance business. In Section \ref{sec:6}, we present some properties of the optimal results and carry on two numerical examples to show the impact of model parameters. Finally, we conclude the paper in Section \ref{sec:7}. 

\section{Model and problem formulation}\label{sec:2}
In this paper, we consider an insurance portfolio which has two dependent classes of business such as motor insurance, health insurance, life insurance, and so on. Let   $(\Omega, \mF,  \pr)$   be a probability space with filtration $ \{ \mF_t \}_{t\geq0}$ containing all objects defined in the following. Let $\{X_i, i\geq 1\}$ be the claim size random variables for the first class  and $\{Y_i, i\geq 1\}$ be the claim size random variables for the second class. Let $X$ and $Y$ be generic random variables which have the same distributions as $X_i$ and $Y_i$, respectively. Then, we assume that  $F_X(x)$ is the common cumulative distribution of  $X_i$ with $F_X(x)=0$ for $x\leq 0$ and $0<F_X(x)\leq1$ for $x>0$, and   $F_Y(y)$  is the common cumulative distribution of $Y_i$ with $F_Y(y)=0$ for $y\leq 0$,  $0<F_Y(y)\leq1$ for $y>0$.
 Then, the aggregate claims up to time $t$ in the two lines of business are denoted by a bivariate compound random process 
$\left(\sum_{i=1}^{M_1(t)}X_i,\sum_{i=1}^{M_2(t)}Y_i\right)$, 
where $M_l(t)=N_l(t)+N(t)$ $(l=1,2)$ is the claim number process for class $l$.  Assume that  $N_1(t)$, $N_2(t)$ and $N(t)$ are three independent Poisson processes with parameters $\lambda_1$, $\lambda_2$ and $\lambda$, respectively;  $X_i$ and $Y_i$ are independent claim size random variables, and that they are independent of $N_1(t)$, $N_2(t)$ and $N(t)$. In this case, the dependence of the two lines of business is  due to a common shock  governed by the counting process $N(t)$. Thus, the total reserve of the insurance portfolio up to time $t$ is given by 
\begin{equation}\label{eq:pro1}
R_t = u + c\,t - \sum_{i=1}^{M_1(t)} X_i-\sum_{i=1}^{M_2(t)} Y_i,
\end{equation}
in which $u  \ge 0$ is the initial surplus and $c$ is the premium rate. 

In this paper, we suppose that the insurer can reinsure its claim with per-loss reinsurance via a continuously payable premium.  Let $\mathcal{H}_{1}=\{\mathcal{H}_{1t}\}_{t\geq0}$ and $\mathcal{H}_{2}=\{\mathcal{H}_{2t}\}_{t\geq0}$  denote the retained claim  at time $t$ for the classes $1$ and $2$, respectively. % thus reinsurance indemnifies the insurer with the amount $X_i-\mathcal{H}_1(X_i)$ for class $1$ and $Y_i-\mathcal{H}_2(Y_i)$ for class $2$.  
 Thus, the surplus of the insurer at time t under the retention  strategy $\mathcal{H}=\{\mathcal{H}_{t}\}_{t\geq0}=(\mathcal{H}_{1},\mathcal{H}_{2})$, denoted by $U=\{U_t\}_{ t  \geq0}$, is given by  
 \begin{equation}\label{eq:Udyn0}
U_t = \left( c-\delta(\mathcal{H}_{t} )\right)t-  \sum_{i=1}^{M_1(t)} \mathcal{H}_{1t}(X_i)-\sum_{i=1}^{M_2(t)} \mathcal{H}_{2t}(Y_i),
\end{equation}
 where  $\delta(\mathcal{H}_{t})$ is the reinsurance premium rate at time $t$ paid to the reinsurer corresponding to the retention policy $\mathcal{H}$. 

We solve the optimization problem by approximating the jump process in \eqref{eq:Udyn0} with a diffusion model.  For convenience, denote the notations $\beta_i$ and $\gamma_i$  $(i=1,2)$ as  follows:
\begin{equation}\label{eq:beta}\beta_1=(\la_1+\la)\E\left(\mathcal{H}_{1}(X_i)\right)
,~~~~ \beta_2=(\la_2+\la)\E\left(\mathcal{H}_{2}(Y_i)\right),\end{equation}
\begin{equation}\label{eq:gamma}\gamma_1=\sqrt{(\la_1+\la)\E\left(\mathcal{H}_{1}(X_i)\right)^2}
,~~~\gamma_2=\sqrt{(\la_2+\la)\E\left(\mathcal{H}_{2}(Y_i)\right)^2}.\end{equation} 
%According to  %Grandell \cite{G91},   Wang and Yuen \cite{WY05},     Bai et al. \cite{BCZ13},  and Liang and Yuen \cite{LY16}, 
Let $\mu_1=\E\mathcal{H}_{1}(X_i)$ and $\mu_2=\E\mathcal{H}_{2}(Y_i).$ Then according to  the proof of Theorem $2.1$ in Bai et al. \cite{BCZ13},   the compound Poisson process $$C_t= \sum_{i=1}^{M_1(t)} \mathcal{H}_{1t}(X_i)+\sum_{i=1}^{M_2(t)} \mathcal{H}_{2t}(Y_i)$$ 
 %is also a compound Poisson process with  parameter $\widetilde{\la}=\la+\la_1+\la_2$, and that the distribution of the transformed claim size random variable $Z$ is given by $$F_Z(z)=\frac{\la_1}{\widetilde{\la}}F_X(z)+\frac{\la_2}{\widetilde{\la}}F_Y(z)+\frac{\la}{\widetilde{\la}}F_{X+Y}(z). $$
 can be  approximated by  the following Brownian motion risk model  
 $$\hat{C}_t= (\beta_1+\beta_2)\, t-\gamma_1dB_{1t}-\gamma_2d B_{2t},$$%= (\beta_1+\beta_2)\, t-\sqrt{\gamma_1^2+\gamma^2_2+2 \rho\gamma_1\gamma_2}\,B_t,$$    
where $B_1=\{B_{1t}\}_{t\geq0}$  and $B_2=\{B_{2t}\}_{t\geq0}$ are two correlated standard Brownian motions with   correlation coefficient 
\begin{equation}\label{eq:rho}\rho=\frac{\la}{\gamma_1\gamma_2}\mu_1\mu_2.\end{equation}
%with $$\mu_1=\E\mathcal{H}_{1}(X_i), ~~~~~~\mu_2=\E\mathcal{H}_{2}(Y_i).$$  %Here,  $ B=\{B_t\}_{t \ge 0}$ is a standard Brownian motion.% which is independent of $W_t$. 

 Furthermore, we assume that the insurer can invest  its surplus in a risky asset  (stock or mutual fund)  and  a  risk-free asset (bond or bank account) with interest rate $r>0$. Specifically,  the process of the  risky asset  follows a geometric Brownian motion given by $$dS_t=\mu S_t\, dt +\sigma S_t\,dW_t,$$
where $\mu>r$ and $\sigma>0$,  and  $W=\{W_t\}_{t\geq0}$ is a standard Brownian motion which is independent of $B_1$ and $B_2$. Let $\pi=\{\pi_t\}_{t\geq0}$  denote the amount invested in the risky asset at time $t\geq0$,  and then the rest of the surplus $(U_t - \pi _t)$ is invested in the risk-free asset. 
%An investment strategy $\pi=\{\pi_t\}_{t\geq0}$ is admissible if it is adapted to the filtration $\{\mF_t\}_{t\geq0}$, and satisfies $\int_0^t\pi_s^2ds<\infty$ with probability one for all $t\geq0$.
%Denote the set of admissible  strategies $(\mathcal{H}_1,\mathcal{H}_2,\pi)$ by $\mA$. 
Thus, for a chosen combination  of  controls $\nu=\{\nu_t\}_{t\geq0}=(\mathcal{H}_1,\mathcal{H}_2,\pi)$, we obtain that the diffusion  approximation surplus process for risk process \eqref{eq:Udyn0} has the following dynamics
%\begin{equation}\label{eq:Udyn}dU_t = \left[ r U_t +(\mu-r)\pi_t+ c-\delta(\mathcal{H}_{t} )\right]dt+\sigma\pi_t\, dW_t -  d\sum_{i=1}^{M_1(t)} \mathcal{H}_{1t}(X_i)-d\sum_{i=1}^{M_2(t)} \mathcal{H}_{2t}(Y_i),\end{equation}where $\delta(\mathcal{H}_{t})$ is  the reinsurance premium rate at time $t$ paid to the reinsurer corresponding to the retention policy $\mathcal{H}_{t}=(\mathcal{H}_{1t},\mathcal{H}_{2t})$. 
%Therefore, the resulting process $\hat{U}=\{\hat{U}_t\}_{ t  \geq0}$ evolves according to the dynamics 
$$
d\hat{U}_t = \left[ r \hat{U}_t +(\mu-r)\pi_t+ c-\delta(\mathcal{H}_{t}) -\beta_1-\beta_2\right]\, dt
+\sigma\pi_t\, dW_t +\gamma_1dB_{1t}+\gamma_2dB_{2t},
$$
or equivalently
 \begin{equation}\label{eq:Uhat}
d\hat{U}_t = \left[ r \hat{U}_t +(\mu-r)\pi_t+ c-\delta(\mathcal{H}_{t}) -\beta_1-\beta_2\right]\, dt
+\sigma\pi_t\, dW_t +\sqrt{\gamma_1^2+\gamma^2_2+2\la\mu_1\mu_2}\, dB_t
\end{equation}
with initial surplus $\hat{U}_0=u$. Here, $B=\{B_t\}_{t\geq0}$ is  a standard Brownian motion which is dependent of $W$.  For  simplicity, we denote \begin{equation}\label{eq:a}a_1=(\la+\la_1)\E X,~~~~~ a_2= (\la+\la_2)\E Y. \end{equation} Then, to avoid mathematical trivialities, we assume that
$
a_1+a_2 < c < \delta(\textbf{0})$, 
   which implies that  the insurer's premium income is greater than the expected value of the claims but less than the premium for full reinsurance.  The notation ``\textbf{0}" denotes the zero vector. Let $\kap$ denote the positive difference 
with  $
\kap =  \delta(\textbf{0}) - c.$ Note that   if the value of the surplus is greater than or equal to
\begin{equation}\label{eq:us}
u_s = \frac{\kappa}{r},
\end{equation}
 then the insurer can buy full reinsurance and invest all the surplus in the risk-free asset to earn interest rate $r$. Then the interest earned  can completely cover the shortfall between premiums and expensive reinsurance purchase, and thus the surplus will never drop below its current value.  For this reason, we call $u_s$ the {\it safe level}. 
 
  In the following, we propose the definition  of admissible strategies.
\begin{definition}
 A strategy $\nu=(\mathcal{H}_1,\mathcal{H}_2,\pi)$ is said to be admissible if the following conditions are satisfied:
 \begin{itemize}
 \setlength{\itemsep}{0pt} 
 \item[(i)]  it is adapted to the filtration $\{\mF_t\}_{t\geq0}$; 
\item [(ii)]  $\mathcal{H}_{it}$  $(i = 1, 2)$ is a  function of the possible claim size $Z=z$ at time $t$, %that is, we can write  $\mathcal{H}_{it}=\mathcal{H}_{it}(x)$ 
and satisfies $0\leq\mathcal{H}_{it}(z)\leq z$, for all $t\geq0$
and $ z\geq0$; 
\item[(iii)]   it satisfies $\int_0^t\pi_s^2ds<\infty$ with probability one for all $t\geq0$.
\end{itemize}
%and satisfies the net profit condition;

\leftline{The set of all admissible strategies is denoted by $\mathcal {A}$.}
\end{definition}

Define the maximum surplus process $M = \{ M_t \}_{t \ge 0}$ by
\begin{equation}\label{eq:M}
M_t = \max \left\{ \sup \limits_{0 \le s \le t} \hU_s, \, M_{0} \right\},
\end{equation}
with $M_0 = m \ge u$.  We allow the surplus process to have a financial past, as embodied by the term $M_{0}$ in \eqref{eq:M}.  Drawdown is the time when the value of the surplus process reaches $\alp \in [0,1)$ times its maximum value, that is, at the hitting time $\tau_\alp$ given by
\begin{equation}\label{eq:tau_alp}
\tau_\alp = \inf\{t \ge 0:\ \hU_t \le \alp M_t \}.
\end{equation}
If $\alp = 0$, then drawdown is the same as ruin for the ruin level $0$.

The minimum probability of drawdown $\phi$ is defined by
\begin{equation}\label{eq:phi}
\phi(u, m) = \inf \limits_{\nu \in \mA} \pr^{u,m}(\tau_\alp < \infty) = \inf \limits_{\nu \in \mA} \E^{u,m} \big({\bf 1}_{\{\tau_\alp < \infty\}} \big),
\end{equation}
in which $\pr^{u,m}$ and $\E^{u,m}$ denote the probability and expectation, respectively, conditional on $\hU_0 = u$ and $M_0 = m$.  
  Note that, if $u \le \alp m$, then $\phi(u, m) = 1$, and if $u \ge u_s$ and $u > \alp m$, then $\phi(u, m) = 0$.  It remains for us to determine the minimum probability of drawdown $\phi$ on the domain
\begin{equation} \label{eq:mO}
\mO = \big\{(u, m) \in (\R^+)^2: \alp m \le u \le \min(m, u_s), \alp m<u_s \big\}.\footnote{Note that if $\alpha m=u_s$, then technically drawdown has occurred, but the insurer could keep its surplus at $\alpha m=u_s$ thereafter by purchasing full reinsurance. Therefore, we avoid this ambiguous case by assuming $\alpha m<u_s$ throughout. }
\end{equation}
\section{Verification theorem} \label{sec:3}
In this section,  we prove a verification theorem, which  is essential in solving the associated stochastic control problem.

 For a given admissible strategy $\nu=(\mathcal{H}_1,\mathcal{H}_2,\pi)$, define  the 
  differential operator $\mA^\nu$ on appropriately differentiable functions by
\begin{equation}\label{eq:op}
\mA^{\nu} h(u,m) = \left[ r u +(\mu-r)\pi+ c-\delta(\mathcal{H}) -\beta_1-\beta_2 \right] h_u + \frac{1}{2}\left(\sigma^2\pi^2+\gamma_1^2+\gamma^2_2+2 \la\mu_1\mu_2\right) h_{uu}. 
 \end{equation}
 Note that the right side of $\eqref{eq:op}$ depends on the maximum surplus $m$ only via $h$'s dependence on m. The verification  theorem is presented as follows:
\begin{theorem}\label{thm:verif}
$($Verification Theorem$)$ Suppose $h: \mO \to \R^+$ is a bounded, continuous function, which satisfies the following conditions:

\noindent
\begin{enumerate}\setlength{\itemsep}{0pt} \setlength{\parsep}{0pt}
\setlength{\parskip}{0pt}
\item[$(i)$] $h(\cdot,m) \in \mC^2((\alp m, \min(m, u_s)))$ is a non-increasing, convex function with bounded first derivative,
\item[$(ii)$] $h(u, \cdot)$ is continuously differentiable, except possibly at $u_s$,
\item[$(iii)$] $h_m(m, m) \ge 0$ if $m < u_s$,
\item[$(iv)$] $h(\alp m, m) = 1$,
\item[$(v)$] $h(u_s, m) = 0$ if $m \ge u_s$,
\item[$(vi)$] $\mA^{\nu}h \ge 0$ for all $\nu \in \mA$.
\end{enumerate}

\noindent
Then, $h \le  \phi$ on $\mO$.

Furthermore, suppose that the function $h$ satisfies all the above conditions in such a way that conditions $(iii)$ and $(vi)$ hold with equality for some admissible strategy $\nu^*=(\mathcal{H}_1^*,\mathcal{H}_2^*,\pi^*)$ defined in feedback form via $\mathcal{H}_{1t}^*=\mathcal{H}_1^*(\hat{U}_t,M_t,X)$, $\mathcal{H}_{2t}^*=\mathcal{H}_2^*(\hat{U}_t,M_t,Y)$ and $\pi_t^*=\pi^*(\hat{U}_t,M_t)$, in which we slightly abuse notation.\footnote{If $m \ge u_s$, then condition $(iii)$ is moot, and we only require equality in condition $(vi)$.} Then, $h = \phi$ on $\mO$, and $\nu^*$ is the optimal reinsurance-investment policy.  
\end{theorem}
\begin{proof} The proof of the verification theorem relies on Lemma 3.1. %whose proof is identical to the one of Appendix A in Han et al. \cite{HLY19}; also see  Lemma $3.1$ in Luo et al. \cite{LWZ19}. 
Then it   can be derived directly  from the corresponding proof given in %in Chen et al. \cite{CLLL15}, Angoshtari et al. \cite{ ABY16b} and 
Han et al. \cite{HLY19}.   We omit the details here. \end{proof}
\begin{lemma}\label{lem:3.1}  Let $\hat{U}_t$ be given in \eqref{eq:Uhat} and $b \in (\alpha m, u_s)$. Define  
$\tau_b =\inf\{t>0: \hU_t\geq b\}$,   $\tau_n=\inf\{t>0: \int_{0}^{t}\pi^2_sds\geq n\},$ $\tau_{\alpha b}=\tau_{\alpha}\wedge\tau_{b}$, and  $\tau=\tau_{\alpha}\wedge\tau_{b}\wedge\tau_{n}.$ 
Then, for any admissible reinsurance-investment  strategy $\nu=\{\mH_{1t},\mH_{2t},\pi_t\}_{t\geq0}$, we have 
$P(\tau_{\alpha b}<\infty) = 1$.
\end{lemma}
\begin{proof} Define the process $Y=\{Y_t\}_{t\geq0}$ by  $Y_t=e^{-\vartheta\hU_t}$. In the expression for $Y_t$, $\vartheta$ is any positive constant to be chosen later in the proof. 
By applying It\^o's  formula to $Y$, we obtain 
\begin{equation}\label{eq:dY}\begin{array}{ll}dY_t=-\vartheta e^{-\vartheta\hU_t}d\hU_t+\frac{1}{2}\vartheta^2 e^{-\vartheta\hU_t }d <\hU,\hU>_t\\[3mm]~~~~=\vartheta e^{-\vartheta\hU_t}\left(\mM_tdt-\sigma\pi_t \,dW_t-\sqrt{\gamma^2_1+\gamma^2_2+2\la\mu_1\mu_2}\,dB_t\right),\end{array}\end{equation}
in which \begin{equation}\label{eq:mM_t}\mM_t=-r\hU_t-c-(\mu-r)\pi_t+\delta(\mH_t)+\beta_1+\beta_2+\frac{1}{2}\vartheta(\sigma^2\pi_t^2+\gamma^2_1+\gamma^2_2+2\la\mu_1\mu_2).\end{equation}
%Assume $\hU_0=u\in(\alpha m,b)$; otherwise, $\tau_{\alpha b}=0$.
 Then according to \eqref{eq:dY}, we have  $$e^{-\vartheta\hU_{\tau\wedge t }}-e^{-\vartheta u}=\int_{0}^{\tau\wedge t}\vartheta e^{-\vartheta\hU_s}\mM_s\,ds-\int_{0}^{\tau\wedge t}\vartheta e^{-\vartheta\hU_s}\sigma \pi_s\,dW_s-\int_{0}^{\tau\wedge t}\vartheta e^{-\vartheta\hU_s}\sqrt{\gamma^2_1+\gamma^2_2+2\la\mu_1\mu_2}\,dB_s.$$
%It follows from the definition of $\tau_n$ that the second integral's expectation equals zero. 
 The integrands of the second and third integrals are bounded, so the integrals' expectation equal zero. Thus, if we take the  expectation on  both sides, we get
$$\E^{u,m} (e^{-\vartheta\hU_{\tau\wedge t}})-e^{-\vartheta u}=\E^{u,m} \int_{0}^{\tau\wedge t}\vartheta e^{-\vartheta\hU_s}\mM_s\,ds.$$
Here, we  first assume that there exists large enough $\vartheta$ such that $\mM_t\geq\frac{\kappa-rb}{2}>0$ for all $0\leq t\leq\tau_{\alpha b}$ with probability $1$. We  give the  specific values of $\vartheta$   in  Lemmas \ref{lem:4.3} and   \ref{lem:4.4} under two different premium principles,  respectively. %which depends on the value of $b$. 
 Note that $e^{-\vartheta b}\leq e^{-\vartheta\hU_{\tau_{\alpha b}\wedge t}}\leq e^{-\vartheta\alpha m}$  for all $t\geq 0$ with probability 1. Then because $e^{-\vartheta\hU_{\tau\wedge t}}$  is bounded by definition, with the above assumption, it follows from the dominated convergence theorem
that 
$$\begin{array}{ll}e^{-\vartheta\alpha m}-e^{-\vartheta u}\geq\E^{u,m} (e^{-\vartheta \hU_{\tau _{\alpha b}\wedge t}})-e^{-\vartheta u}=\E^{u,m} \displaystyle\int_{0}^{\tau_{\alpha b}\wedge t}\vartheta e^{-\vartheta\hU_s}\mM_s\,ds\\[3mm]\geq \vartheta e^{-\vartheta b}\displaystyle\frac{\kappa-rb}{2}\E^{u,m}\int_{0}^{\tau_{\alpha b}\wedge t}{ 1}\,ds
\\[3mm]=\vartheta e^{-\vartheta b}\displaystyle\frac{\kappa-rb}{2}\left(\E^{u,m} \displaystyle\int_{0}^{\tau_{\alpha b}}1 \,ds \cdot{\bf 1}_{\{\tau_{\alpha b}\leq t\}} +\E^{u,m} \displaystyle\int_{0}^{t}1\,ds \cdot{\bf 1}_{\{ \tau_{\alpha b}>t\}} \right)\\[3mm]\geq \vartheta e^{-\vartheta b}\displaystyle\frac{\kappa-rb}{2} t\, \pr^{u,m}(\tau_{\alpha b}>t).\end{array}$$
Letting $t\rightarrow \infty$ in the last expression, we deduce that  $\lim_{t \to \infty}\pr^{u,m}(\tau_{\alpha b}>t)=0$, or equivalently, $\lim_{t \to \infty}\pr^{u,m}(\tau_{\alpha b}<\infty)=1$. We complete the proof.
\end{proof}
%\begin{remark}\label{rem:3.1} We point out  that Lemma \ref{lem:3.1} is indispensable for proving the verification theorem.  According to  the discussion in the end of  Section \ref{sec:2},  both the drift and the volatility of the controlled surplus process in \eqref{eq:Uhat} approach $0$  as the surplus approaches $u_s$. Thus, our intuition tells us that the safe level might not be reachable, and  Bayraktar and Zhang {\rm \cite{BZ15}},  Angoshtari et al.\ {\rm \cite{ABY16a}}  and Han et al .\cite{HLY19} showed that    optimally controlled surplus does not reach the safe level  with positive probability in finite time, which confirms our intuition. 
% \end {remark}

Recall the definition of the domain $\mO$ in \eqref{eq:mO}.  In the case of $m\geq u_s$,  we must have either $\hat{U}_t<u_s$ almost surely
 for all $t\geq0$, or $\hat{U}_t=u_s$ for some $t>0$.
  Thus, the maximum level of surplus $M_t$ does not increase and avoiding drawdown is equivalent to avoiding ruin with a (fixed) ruin level of $\alpha m$.
  However, in the  case of  $m\leq u_s$, $M_t$ can be larger than $m$, i.e.,
the drawdown level becomes a non-decreasing process.
Therefore, according to the verification theorem, when  $\alpha m\leq u\leq u_s\leq m$,  we try to find a smooth, decreasing, convex solution $\xi$  solves  the following differential equation 
\begin{equation}\label{eq:b_pro}
\min \limits_{\nu \in \mA}\left\{\displaystyle \left[ r u +(\mu-r)\pi_t+ c-\delta(\mathcal{H}_{t}) -\beta_1-\beta_2 \right] \xi_u + \frac{1}{2}\left(\sigma^2\pi^2+\gamma_1^2+\gamma^2_2+2 \la\mu_1\mu_2\right) \xi_{uu}\right\}=0,
\end{equation}
with   the  boundary conditions
$$
\xi(\alpha m,m)=1,~\xi(u_s,m)=0.$$
When $\alpha m\leq u\leq m< u_s$, the candidate solution $\xi$ should  still solve the equation \eqref{eq:b_pro} but satisfy the boundary conditions 
$$
 \left\{\begin{array}{ll}
\xi(\alpha m,m)=1,~~~~~ \xi(u_s,u_s)=0,\\[3mm]
\xi_m(m,m)=0.
\end{array}\right.
$$
\section{Optimal  strategy and  value function}\label{sec:4}

In this section, we compute the optimal reinsurance and investment strategies  to minimize the probability of drawdown on region $\mO$  for  the risk model \eqref{eq:Uhat}.  Combining with the verification theorem in Section 3, the optimal results   are derived not only for the expected value   principle but also for the variance premium principle. 

\subsection{The expected value   principle}\label{sec:4.1}
When the reinsurance premium is calculated by  the expected value principle,  recall  from the definition of
$\beta_i$  and $a_i$ $(i=1,2)$ in  \eqref{eq:beta} and \eqref{eq:a},  the reinsurance premium rate at time $t$ is given by 
\begin{equation}\begin{array}{ll}\label{eq:ex-pre}\delta(\mH_t)%=(1+\eta_1)(\la+\la_1)\E (X_i-\mH_1(X_i))+(1+\eta_2)(\la+\la_2)\E (Y(i)-\mH_2(Y_i))\\[3mm]~~~~~~~~
=(1+\eta_1)\left(a_1-\beta_1\right)+(1+\eta_2)\left(a_2-\beta_2\right),\end{array}\end{equation}
where $\eta_1$ and $\eta_2$ are reinsurer's safety loading of the two classes of the insurance business, respectively. Without loss of generality, we assume that $\eta_1\geq\eta_2$ since   the result for  the
case of $\eta_1<\eta_2 $ can be obtained  along the same lines. 

We substitute $\delta(\mH_t)$ above into \eqref{eq:b_pro} and  define  a related  function  
\begin{equation}\label{eq:g1}\begin{array}{lll}g_1(u,\mH_1,\mH_2,\pi) =\big(ru-\kappa+(\mu-r)\pi +\eta_1\beta_1+\eta_2\beta_2\big)\xi_u+\frac{1}{2}(\sigma^2\pi^2+\gamma^2_1+\gamma^2_2+2\la\mu_1\mu_2)\xi_{uu}.
\end{array}
\end{equation}
By using the cumulative distribution functions of $X$ and $Y$, respectively,  we  rewrite $g_1$ as follows:
\begin{equation}\label{eq:ng1}\begin{array}{lll}g_1(u,\mH_1,\mH_2,\pi) =\big(ru-\kappa+
(\mu-r)\pi+\eta_2(\la+\la_2)\E \mH_2\big)\xi_u+\displaystyle\frac{1}{2}\big(\sigma^2\pi^2+(\la+\la_2)\E\mH^2_2\big)\xi_{uu}

\\[5mm]~~~~~~~~~~~~~~~~~~~~~~~ +\displaystyle\int_0^\infty \left\{(\la+\la_1)\eta_1 \mH_1(x)\xi_{u}+\Big(\frac{1}{2}(\la+\la_1)\mH_1^2(x)+\la\mH_1(x)\E\mH_2\Big)\xi_{uu}\right\}dF_X(x)
\\[5mm]~~~~~~~~~~~~~~~~~~~~~=\big(ru-\kappa+
(\mu-r)\pi+\eta_1(\la+\la_1)\E \mH_1\big)\xi_u+\displaystyle\frac{1}{2}\big(\sigma^2\pi^2+(\la+\la_1)\E\mH^2_1\big)\xi_{uu}\\[5mm]~~~~~~~~~~~~~~~~~~~~~~~ +\displaystyle\int_0^\infty \left\{(\la+\la_2)\eta_2 \mH_2(x)\xi_{u}+\Big(\frac{1}{2}(\la+\la_2)\mH_2^2(x)+\la\mH_2(x)\E\mH_1\Big)\xi_{uu}\right\}dF_Y(x).
\end{array}
\end{equation}
From this integral representation of $g_1$, we deduce that we can minimize $g_1$ by minimizing the integrand $x$-by-$x$, subject to $0 \le \mH_i(x) \le x$ $(i=1,2)$.
As a function of $\mH_1$,  $\mH_2$ and $\pi$,
 the integrand is a parabola, so it is minimized  by\footnote{ We slightly abuse the notions of $\mH^*_i (i=1, 2)$ by adding an argument $u$. }

$$\mH^*_1(u,x)=d^*_1(u) \wedge x, ~~~~~~\mH^*_2(u, y)=d^*_2(u) \wedge y$$with 
\begin{equation}\label{eq:d1*}d^*_1(u)=-\eta_1\frac{\xi_u}{\xi_{uu}}-\frac{\la}{\la+\la_1}\E\mH^*_2,\end{equation}
 \begin{equation}\label{eq:d2*}d^*_2(u)=-\eta_2\frac{\xi_u}{\xi_{uu}}-\frac{\la}{\la+\la_2}\E\mH^*_1,\end{equation}   % Since the price process of risky asset is process is independent of   the claim process, we have
 and  \begin{equation}\label{eq:pi*}\pi^*(u)=-\displaystyle\frac{u-r}{\sigma^2}\frac{\xi_u}{\xi_{uu}}.\end{equation}
 To simplify further analysis, we define the following functions
$$h_X(d) =  \E [X_i\wedge d] = \int_0^{{d}} S_X(x) dx,$$
$$h_Y(d) =  \E [Y_i\wedge d] = \int_0^{{d}} S_Y(x) dx,$$
$$H_X(d) =  \E[X_i\wedge d] ^2= 2 \int_0^{d} x S_X(x) dx,$$
$$H_Y(d) =  \E[Y_i\wedge d] ^2= 2 \int_0^{d} x S_Y(x) dx, $$
where $S_X(x)=1-F_X(x)$ and  $S_Y(x)=1-F_Y(x)$. According to  the expressions given in \eqref{eq:d1*}-\eqref{eq:pi*}, it is straightforward to see that the  point   $(d^*_1,d^*_2,\pi^*)\footnote{These valuables depend on the argument $u$, but for simplicity of notation, we omit it in some expressions.}$  should  satisfy the following equations \begin{equation}\label{eq:d1}(\la+\la_1)\eta_1 \xi_u+[(\la+\la_1)d^*_1+\la h_Y(d^*_2)]\xi_{uu}=0,\end{equation}
\begin{equation}\label{eq:d2}(\la+\la_2)\eta_2 \xi_u+[(\la+\la_2)d^*_2+\la h_X(d^*_1)]\xi_{uu}=0,\end{equation}
\begin{equation}\label{eq:pi}(\mu-r)\xi_u+\sig^2\pi^* \xi_{uu}=0.\end{equation}

Note that the Hessen matrix of  ${g_1}$  at point $(\mH^*_1,\mH^*_2,\pi^*)$ can be written as  
\begin{equation}\label{eq:matrixD}\textbf{D}=\xi_{uu}\left(\begin{array}{ccccc}
&(\la+\la_1)S_X(d^*_1)& \la S_X(d^*_1)S_Y(d^*_2) &0\\[3mm]
&\la S_X(d^*_1)S_Y(d^*_2)&(\la+\la_2)S_Y(d^*_2) &0\\[3mm]
&0&0 & \sigma^2
 \end{array}\right),\end{equation}
which is clearly  positive definite.  Therefore, if we can find a convex candidate solution $\xi$ and a point  $(d^*_1,d^*_2,\pi^*)$  such that \eqref{eq:d1}-\eqref{eq:pi} hold, then the point $(\mH^*_1,\mH^*_2,\pi^*)$ is the extreme minimum point of $g_1$. In particular, we point  out that  $d^*_1$, $d^*_2$ and $\pi^*$ have the following relationship
%In order to find the optimal policy and the value function, we need to study the solutions to \eqref{eq:d1} and \eqref{eq:d2} in following context. Once we obtain the values of   $d^*_1$ and $d^*_2$, then $\pi^*$ can be derived by the following equations
\begin{equation}\label{eq:xiuu/xiu}\begin{array}{ll}\displaystyle-\frac{\xi_{uu}}{\xi_{u}}=\frac{(\la+\la_1)\eta_1}{(\la+\la_1)d^*_1+\la h_Y(d^*_2)}
\\[7mm]~~~~=\displaystyle\frac{(\la+\la_2)\eta_2}{(\la+\la_2)d^*_2+\la h_X(d^*_1)}
=\displaystyle\frac{u-r}{\sigma^2\pi^*}.\end{array}\end{equation}
Inspired by Bai et. al \cite{BCZ13}, we define the following auxiliary functions \begin{equation}\label{eq:lxly}l_X(x)=\eta_2 x-\frac{\la\eta_1}{(\la+\la_2)}h_{X}(x),~~~~~l_Y(x)=\eta_1 x-\frac{\la\eta_2}{(\la+\la_1)}h_{Y}(x),\end{equation}
and \begin{equation}\label{eq:k} k(u,x)=(ru-\kappa)+(\la+\la_1)\eta_1\Big(h_X(x)-\frac{H_X(x)}{2x}\Big)+\frac{(\mu-r)^2x}{2\eta_1\sigma^2}.\end{equation}
It is not difficult to find that  $l_X(d^*_1)=l_Y(d^*_2),$ and $l_X(0)=l_Y(0)=0.$ Because of $\eta_1\geq \eta_2$,  it follows that $l_Y(x)$ is  a strictly  increasing function with respect to $x\geq0$,  and thus the inverse function  $l^{-1}_Y(y)$ exists and   strictly increases with respect to $y\geq 0$  and satisfies  $l^{-1}_Y(0)=0$. Besides, it is  clear that  the function $k(u,x)$ is also a  strictly increasing function with respect to  $x$, and hence for any fixed $u$, the inverse function $k^{-1}(u,y)$ exists. It turns out that  $k(u,0)<0$ and $\lim_{x \to \infty} k(u,x) =+\infty$ hold for all $ \alp m \le u \le \min(m, u_s)$.  Correspondingly, we define two zeros associated with $l_X(x)$ and $k(u,x)$ respectively as follows
\begin{equation}\label{eq:a_l} a_l=\sup\left\{x\geq0: l_X(x)=0\right\},\end{equation}
and 
\begin{equation}\label{eq:a_k}a_k(u)=k^{-1}(u,0).\end{equation}
Note that  both $a_l$ and $a_k$ are nonnegative. It is easy to check that $a_k(u)$ is a strictly decreasing function with  respect to $u$ and satisfies  the boundary condition $a_k(u_s)=0$. In particular, there exists a $\widetilde{u}\leq u_s$ such that $a_k(u)\geq a_l$ when $u\leq\widetilde{u}$ and  $a_k(u)< a_l$ when $u>\widetilde{u}$.  Since we are only interested in the domain $\mO$, without loss of generality, we assume that  $ \widetilde{u}\geq \alpha m$. Here, we point out that if   $ \widetilde{u}<\alpha m$, we have $a_k(u)< a_l$ for all $u\in[\alpha m, \min(m,u_s)]$, and thus  the optimal results   can be directly obtained  from  the case of $\widetilde{u}< u\leq \min(m,u_s)$ by replacing $\widetilde{u}$ with $\alpha m$.

In order to find the value function and the optimal reinsurance-investment strategy, we study the solutions to \eqref{eq:d1}-\eqref{eq:pi} in the following lemmas.
% $$k(x)=(ru-\kap)+(\la+\la_1)\eta\left(h_X(x)-\frac{H_X(x)}{2x}\right)+\frac{(\mu-r)^2x}{2\eta\sigma^2}$$
\begin{lemma}\label{lem:4.1}
Suppose that  $\alpha m\leq u\leq \min(m,\widetilde{u})$.  Then  there  exists a   point $(d^*_1,d^*_2,\pi^*) $ which  uniquely solves  the system \eqref{eq:d1}-\eqref{eq:pi}. 
\end{lemma}

\noindent
{\bf Proof:} Because of  $l_X(d^*_1)=l_Y(d^*_2)$ and  the inverse function  $l^{-1}_{Y}$ exists,  it follows that  \begin{equation}\label{eq:d1=d2}d^*_2=l^{-1}_{Y}(l_X(d^*_1)).\end{equation}  
Putting \eqref{eq:xiuu/xiu} and \eqref{eq:d1=d2} back into \eqref{eq:g1}  and letting $g_1=0$ yield
\begin{equation}\label{eq:eqG}G(u, d^*_1)\xi_u=0,\end{equation}
where the function $G(u,d)$ is given by
\begin{equation}\label{eq:G}\begin{array}{ll}G(u,d)=( r u -\kappa) +\eta_1(\la+\la_1)h_X(d)+\eta_2(\la+\la_2)h_Y(l^{-1}_Y(l_X(d)))\\[5mm]~~~~~~~~~~~~~-\displaystyle\frac{\eta_1}{2} \frac{ (\la+\la_1)H_X(d)+(\la+\la_2)H_Y\big(l^{-1}_Y(l_X(d))\big)+2\la h_X(d)h_Y\big(l^{-1}_Y(l_X(d))\big)}{d+\frac{\la}{\la+\la_1}h_Y(l^{-1}_Y(l_X(d)))}\\[5mm]~~~~~~~~~~~~~~+\displaystyle\frac{1}{2}\frac{(\mu-r)^2\left(d+\frac{\la}{\la+\la_1}h_Y\big(l^{-1}_Y(l_X(d))\big)\right)}{\eta_1\sigma^2}. 
\end{array}\end{equation}
%In the following, we prove that there exists a unique  positive solution $d^*_1$ to the equation $G(d)=0$.
%Here, we slightly abuse notation by replacing  the arguments $(u,\mH_1,\mH_2,\pi)$  of $g_1$ by $(u, d_1)$.   
Note that for any $d\geq a_l$, 
$$\begin{array}{ll}\displaystyle\frac{\partial G(u, d)}{\partial d}=\displaystyle\frac{\eta_1(\la+\la_1)}{2}\big[(\la+\la_1)H_X(d)+(\la+\la_2)H_Y(l^{-1}_Y(l_X(d)))+2\la h_X(d)h_Y(l^{-1}_Y(l_X(d)))\big]\\[5mm]~~~~~~~~~~~~~~\times\displaystyle\frac{(\la+\la_1)+\la S_Y(l^{-1}_Y(l_X(d)))(l^{-1}_Y(l_X(d)))^\prime l^\prime_X(d)}{[(\la+\la_1)d+\la h_Y(l^{-1}_Y(l_X(d)))]^2}\\[5mm]~~~~~~~~~~~~~~~+\displaystyle\frac{(\mu-r)^2\left[(\la+\la_1)+\la S_Y(l^{-1}_Y(l_X(d)))(l^{-1}_Y(l_X(d)))^\prime l^\prime_X(d)\right]}{2(\la+\la_1)\eta_1\sigma^2}>0,\end{array}$$
which implies that $ G(u, d)$ is a strictly increasing function with respect to $d$ for $d\in[a_l, \infty)$. Besides, since  $l_X(a_l)=l^{-1}_Y(l_X(a_l))=0$,  combining with  the definition of $k(u,x)$ in \eqref{eq:k},   we have
\begin{equation}\label{eq:G(al)}\begin{array}{ll}G(u, a_l)=\displaystyle(ru-\kappa)+(\la+\la_1)\eta_1\Big(h_X(a_l)-\frac{H_X(a_l)}{2a_l}\Big)+\frac{(\mu-r)^2a_l}{2\eta_1\sigma^2}\\[3mm]~~~~~~~~~~~=k(u,a_l).\end{array}\end{equation}
When $\alpha m\leq u\leq\min(m, \widetilde{u})$, one can show that  $a_l\leq a_k(u)$, and thus it follows that $k(u,a_l)<0$. Moreover,  it is not difficult to find that 
$$\lim_{d \to \infty}\displaystyle\frac{1}{2}\frac{(\mu-r)^2\left(d+\frac{\la}{\la+\la_1}h_Y\big(l^{-1}_Y(l_X(d))\big)\right)}{\eta\sigma^2} =+\infty, $$
then we have 
$\lim_{d \to \infty} G(u, d) =+\infty.$ Therefore, $G(u, d)$ has a unique zero at $ d_1^*(u)\in[a_l,\infty)$ for  all $u\in[\alp m,   \min(m,\widetilde{u})]$, from which it follows that the equations given in \eqref{eq:d1}-\eqref{eq:pi} have a unique solution $(d^*_1,d^*_2,\pi^*)$  with 
\begin{equation}\label{eq:d2*}d^*_2(u)=l^{-1}_Y\big(l_X(d^*_1(u))\big),\end{equation} and
$$\pi^*(u)=\frac{(\mu-r)\Big(d^*_1(u)(\la+\la_1)+\la h_Y(l^{-1}_Y\big(l_X(d^*_1(u))\big)\big)\Big)} {\eta_1\sigma^2(\la+\la_1)}.$$ We complete the proof.\qed
%\begin{remark} 

\vspace{3mm}
Note  that  when $\widetilde{u}<u\leq\min(m, u_s)$, we always have $a_l>a_k(u)$, which implies that  $G(u, d)>0$ for $d=a_l$. According to the monotonicity of  $G(u, d)$,  there does not  exist $d^*_1\geq a_l$ and $d^*_2\geq 0$  such that  \eqref{eq:eqG} and \eqref{eq:d2*} hold.   Thus, the extreme minimum point  of $g_1$ in \eqref{eq:g1}  can only be obtained on $\mH^*_1(u,x)=\bar{d}^*_1(u)\wedge x$ with $ \bar{d}^*_1(u)\in[0, a_l)$ and  $\mH^*_2(u,y)=\bar{d}^*_2(u)\wedge y$ with $\bar{d}^*_2(u)=0$. Therefore, we have the following lemma.
%\end{remark}
%We give  the extreme minimum  point $(\mH^*_1,\mH^*_2, \pi^*)$  of $g_1$  in following  Lemma.
\begin{lemma}\label{lem:4.2}
When $\widetilde{u}<u\leq\min(m, u_s)$,  the extreme minimum  point   of $g_1$ is given by $$\mH^*_1(u, x)=a_k(u)\wedge x,~~~~~\mH^*_2(u, y)=0,$$ where $a_k(u)$ is defined in \eqref{eq:a_k} and $$\pi^*(u)= \frac{\mu-r}{\sig^2}\cdot\frac{a_k(u)}{\eta_1}.$$ \end{lemma}

\noindent
{\bf Proof:} When $\widetilde{u}<u\leq\min(m, u_s)$,  putting $\mH^*_2=0$ back into \eqref{eq:ng1} and differentiating $g_1$
 with respect to $\mH_1$ and $\pi$ yield 
 $$\mH^*_1(u,x)=\bar{d}^*_1(u) \wedge x$$
 with \begin{equation}\label{eq:d1*_2}\bar{d}^*_1(u)=-\eta_1\frac{\xi_u}{\xi_{uu}},\end{equation}
and $ \pi^*(u)$ given in the form of \eqref{eq:pi*}.
Along the same lines as in Lemma \ref{lem:4.1}, we finally derive  that   $\bar{d}^*_1$ should satisfy  the equation $k(u,\bar{d}^*_1)=0$. Since $k(u,x)$ is a strictly increasing function with respect to $x$ and the inequalities $k(u,0)<0$  and $k(u, a_l)>0$ hold for all $\widetilde{u}<u\leq\min(m, u_s)$, we have
$\bar{d}^*_1(u)=a_k(u)$, which completes the proof.
\qed
\begin{remark}\label{re:us}
We can see from Lemma \ref{lem:4.2} that   when the value of the surplus increases toward $u_s$, the insurer optimally  buys  full reinsurance and invests all the surplus in the risk-free asset due to the fact that $a_k(u)=0$ for $u=u_s$. As a result, wealth will never decrease and drawdown cannot happen, which   is consistent with the discussion at the end of Section \ref{sec:2}. 
\end{remark}
In the following lemma, as mentioned in Lemma \ref{lem:3.1}, we will prove that   there indeed  exists  large enough $\vartheta$ such that $\mM_t\geq\frac{\kappa-rb}{2}>0$ for all $0\leq t\leq\tau_{\alpha b}$ with probability 1.
\begin{lemma}\label{lem:4.3} Recall   $\mM_t$  of  \eqref{eq:mM_t} and $b \in (\alpha m, u_s)$. Let $\delta(\mH_t)$ be given in \eqref{eq:ex-pre}.    Then, there exists $\vartheta\in(0,\infty)$  given by
$$\vartheta=\left\{\begin{array}{llll}\eta_1\left(\widetilde{d}_1+\frac{\la}{\la+\la_1}h_Y(l^{-1}_Y(l_X(\widetilde{d}_1)))\right)^{-1},~~~~~& \alpha m\leq b \leq\max(\alpha m, 2\widetilde{u}-u_s),\\[4mm]\eta_1\left(a_k\Big(\frac{b}{2}+\frac{\kappa}{2r}\Big)\right)^{-1},~~~~~&\max(\alpha m, 2\widetilde u-u_s)<b\leq u_s \end{array}\right.$$
such that $\mM_t\geq\frac{\kap-rb}{2}.$ Here, $\tilde{d}_1\in[a_l, \infty)$ uniquely solves $G\Big(\frac{b}{2}+\frac{\kappa}{2r},d\Big)=0$
with $G$  defined in \eqref{eq:G}, and $a_k$ is defined in \eqref{eq:a_k}  \end{lemma}
\begin{proof}  Substituting \eqref{eq:ex-pre}  back into the expression of $\mM_t$, it follows that 
 $$\mM_t=\kappa-r\hU_t-(\mu-r)\pi_t-\eta_1\beta_1-\eta_2\beta_2+\frac{1}{2}\vartheta(\sigma^2\pi_t^2+\gamma^2_1+\gamma^2_2+2\la\mu_1\mu_2).$$
%We wish to show that, for $\delta$ large enough, $\mM_t\geq\frac{\kappa-rb}{2}>0$ for all $t\geq0$ with probability 1. 
Since  $\hU_t\leq b<u_s=\kappa/r$ for all $0\leq t\leq\tau_{\alpha b}$, we have $\kap-r\hU_t\geq \kap-rb>0$. We then find $\vartheta>0$ such that 
$$\begin{array}{ll}%\Lambda\left((\la+\la_1)\E(\mH_1^2-2\mH_1X)+(\la+\la_2)\E(\mH_2^2-2\mH_2Y)+2\la\E(\mH_1\mH_2-\mH_1Y-\mH_2X)\right)
-\eta_1\beta_1-\eta_2\beta_2-(\mu-r)\pi+\displaystyle\frac{1}{2}\vartheta(\sigma^2\pi^2+\gamma^2_1+\gamma^2_2+2\la\mu_1\mu_2)\geq-\frac{\kap-rb}{2},\end{array}$$
or equivalently, \begin{equation}\label{eq:ex_ine1}\frac{1}{2}\vartheta\geq\frac{%-\Lambda\left((\la+\la_1)\E(\mH_1^2-2\mH_1X)+(\la+\la_2)\E(\mH_2^2-2\mH_2Y)+2\la\E(\mH_1\mH_2-\mH_1Y-\mH_2X)\right)
\eta_1\beta_1+\eta_2\beta_2+(\mu-r)\pi-\frac{\kap-rb}{2}}{\sigma^2\pi^2+\gamma^2_1+\gamma^2_2+2\la\mu_1\mu_2}\end{equation}
for all admissible strategies $\nu=(\mH_1,\mH_2,\pi)$.
In other words, we wish to show that the right side of  inequality \eqref{eq:ex_ine1} has a finite maximum over such $\nu$. By comparing the right side of the inequality with ${g}_1$ in \eqref{eq:g1}, this problem is essentially to the ones we solved in Lemmas \ref{lem:4.1}-\ref{lem:4.2}. By relying on the proof of those lemmas, we deduce that, when  $ \alpha m\leq b \leq\max(\alpha m, 2\widetilde{u}-u_s)$, the right side of the inequality is maximized by
$$\tilde{\mH}_1(x)=\tilde{d}_1 \wedge x, ~~~~~~\tilde{\mH}_2(y)=l^{-1}_{Y}(l_X(\tilde{d}_1)) \wedge y,$$
 and  $$\tilde{\pi}=\frac{(\mu-r)\Big(\tilde{d}_1(\la+\la_1)+\la h_Y\big(l^{-1}_Y(l_X(\tilde{d}_1))\big)\Big)} {\eta_1\sigma^2(\la+\la_1)},$$
in which $\tilde{d}_1\in[a_l, \infty)$ uniquely solves $$G\Big(\frac{b}{2}+\frac{\kappa}{2r},d\Big)=0.$$
By substituting $(\tilde{\mH}_1,\tilde{\mH}_2,\tilde{\pi})$  into the right side of \eqref{eq:ex_ine1}, we obtain  $$\frac{\vartheta}{2}>\frac{\eta_1}{2}\cdot\big(\tilde{d}_1+\frac{\la}{\la+\la_1}h_Y(l^{-1}_Y(l_X(\tilde{d}_1)))\big)^{-1}.$$Thus, set $\vartheta=\eta_1\big(\tilde{d}_1+\frac{\la}{\la+\la_1}h_Y(l^{-1}_Y(l_X(\tilde{d}_1)))\big)^{-1}\in(0,\infty),$ and for this choice of $\vartheta$, we have $\mM_t\geq\frac{\kap-rb}{2}>0$ for all $0\leq t\leq\tau_{\alpha b}$ with probability 1.
When  $\max(\alpha m, 2\tilde u-u_s)<b\leq u_s$,  the right side of the inequality is maximized by
 $$\tilde{\mH}_1(x)=a_k\big(\frac{b}{2}+\frac{\kappa}{2r}\big)\wedge x,~~~~~\tilde{\mH}_2(y)=0,$$ and $$\tilde{\pi}= \frac{\mu-r}{\sig^2}\cdot\frac{a_k\big(\frac{b}{2}+\frac{\kappa}{2r}\big)}{\eta_1}.$$
By substituting $(\tilde{\mH}_1,\tilde{\mH}_2,\tilde{\pi})$  into the right side of \eqref{eq:ex_ine1}, we obtain  $$\frac{\vartheta}{2}>\frac{\eta_1}{2}\Big(a_k\big(\frac{b}{2}+\frac{\kappa}{2r}\big)\Big)^{-1}.$$Thus, set $\vartheta=\eta_1\Big(a_k\Big(\frac{b}{2}+\frac{\kappa}{2r}\Big)\Big)^{-1}\in(0,\infty),$ and for this choice of $\vartheta$, we have $\mM_t\geq\frac{\kap-rb}{2}>0$ for all $0\leq t\leq\tau_{\alpha b}$ with probability 1. We complete the proof.
\end{proof}
In the following theorem, according to the analysis  given in Lemmas \ref{lem:4.1} and \ref{lem:4.2}, we compute our candidate minimum probability   of drawdown for  both the cases of $m\geq u_s$ and $m<u_s$. Then combining with the verification theorem given in Section \ref{sec:3}, we verify that the resulting functions are indeed the minimum probability of drawdown on $\mO$.

\begin{theorem}\label{th3.1} Suppose that   $\tilde{u}\geq\alpha m$. Let $g_i$  and $f_i$ $(i=1,2)$ be given in Appendix   A.1. If $u_s \leq m$,   for any $u\in[\alpha m, u_s]$, the minimum probability of drawdown for the surplus process \eqref{eq:Uhat} is given by \begin{equation} \label{eq:drawdown1}
 \phi(u,m)=\left\{\begin{array}{lll}
 \displaystyle1-\frac{g_{1}(u,m)}{g_{2}(u_s,m)},& \alpha m\leq u \leq\widetilde{u}, \\[5mm]
 \displaystyle1-\frac{g_{2}(u,m)}{g_{2}(u_s,m)}, & \widetilde u<u\leq u_s\leq m.
 \end{array}\right.
\end{equation}
For $m< u_s$,
\noindent
 (\romannumeral1) if $\widetilde{u} \leq m < u_s$, then for any $u\in[\alpha m, m]$, the minimum probability of drawdown for the surplus process \eqref{eq:Uhat} is given by
\begin{equation}\label{eq:drawdown2}
 \phi(u,m)=\left\{\begin{array}{llll}
  \displaystyle1-k_{2}(m)\cdot\frac{g_{1}(u,m)}{g_{2}(u_s,u_s)}, & \alpha m\leq u\leq \widetilde{u},\\[5mm]
 \displaystyle1-k_{2}(m)\cdot\frac{g_{2}(u,m)}{g_{2}(u_s,u_s)}, & \widetilde{u}< u\leq m<u_s,
 \end{array}\right.
\end{equation}
where
$$
k_{2}(m)=\exp\left\{-\int_m^{u_s}{f_{2}(y)dy}\right\};
$$
\noindent
(\romannumeral2) if $\alpha m <m < \widetilde{u}$, then for any $u\in[\alpha m,m]$, the minimum probability of drawdown for the surplus process \eqref{eq:Uhat} is given by
\begin{equation}\label{eq:drawdown3}
 \phi(u,m)=\displaystyle1-k_{1}(m)\cdot\frac{g_{1}(u,m)}{g_{2}(u_s,u_s)},
\end{equation}
where
$$
k_{1}(m)=\exp\left\{-\left(\int_{m}^{\widetilde{u}}{f_{1}(y)}+\int_{\widetilde{u}}^{u_s}{f_{2}(y)}\right)dy\right\}.
$$
The corresponding  optimal reinsurance-investment strategy, for all $t\geq0$, is given in feedback form by 
\begin{equation}\label{eq:strategy1}\left \{\begin{array}{ll}\mH^*_{1t}=d^*_1(\hat{U}_t)\wedge X,~~~~~~\mH^*_{2t}=d^*_2(\hat{U}_t)\wedge Y,\\[5mm]
\pi^*_t=\displaystyle\frac{\mu-r}{\sigma^2}\left(\frac{d^*_1(\hat{U}_t)}{\eta_1}+\frac{\la h_Y(d^*_2(\hat{U}_t))}{\eta_1(\la+\la_1)}\right)
\end{array}\right. \end{equation}
for all $u\in [\alpha m, \min(m,\widetilde{u})]$, and 
\begin{equation}\label{eq:strategy2}\left \{\begin{array}{ll}\mH^*_{1t}=a_k(\hat{U}_t)\wedge X,~~~~~~
\mH^*_{2t}=0,\\[5mm]
\pi^*_t=\displaystyle\frac{\mu-r}{\sigma^2}\cdot\frac{a_k(\hat{U}_t)}{\eta_1}
\end{array}\right.\end{equation}
for all $u\in (\widetilde{u},\min(m,u_s)]$.
 Here, $d^*_1(u)\geq a_l$ uniquely solves \eqref{eq:G}, $d^*_2(u)=l^{-1}_{Y}(l_X(d^*_1(u)))$ and $a_k(u)$ is defined by \eqref{eq:a_k}.\end{theorem}

\noindent
{\bf Proof:} When $m\geq u_s$,   $\xi$ solves the boundary-value problem
\begin{equation}\label{eq:b-v}\left\{\begin{array}{ll}\displaystyle(\la+\la_1)\eta_1\xi_u+\Big[{(\la+\la_1)d^*_1+\la h_Y(d^*_2)}\Big]\xi_{uu}=0,  &~\alpha m \leq u\leq \widetilde{u}, \\[3mm]

\eta_1\xi_u+a_k(u)\xi_{uu}=0, &~\widetilde{u}<u\leq u_s, \end{array}\right.\end{equation}
with boundary conditions 
$$
\xi(\alpha m,m)=1,~~~~~ \xi(u_s,m)=0.$$
Let $\xi$ equal to the right-hand side of \eqref{eq:drawdown1}, then we can clearly see that $\xi$ solves this boundary-value problem. Furthermore, it is straightforward to show that $\xi$ satisfies  the conditions of  Theorem  \ref{thm:verif}. In particular, Condition $(\romannumeral 3)$ is moot because of  $m\geq u_s$ and  Condition $ (\romannumeral 6)$ holds  with equality.  Thus, the probability of drawdown $\phi$ for $m\geq u_s$ is indeed the expression given in \eqref{eq:drawdown1} associated with the reinsurance-investment strategy given in  \eqref{eq:strategy1} and \eqref{eq:strategy2}.

When $m<u_s$, $\xi$ still solves the differential equations \eqref{eq:b-v} but with  the boundary conditions $$
 \left\{\begin{array}{ll}
\xi(\alpha m,m)=1,~~~~~ \xi(u_s,u_s)=0,\\[3mm]
\xi_m(m,m)=0.
\end{array}
\right.$$
Let $\xi$ equal to the right-hand sides  of \eqref{eq:drawdown2} and   \eqref{eq:drawdown3}. It is clear to  verify that the $\xi$ solves this boundary-value problem.  Furthermore, we can see that $\xi$ satisfies all  the conditions of  Theorem  \ref{thm:verif}; specifically,  $\xi$ satisfies Conditions $(\romannumeral 3)$ and $ (\romannumeral 6)$  with equality.  Thus,  $\xi$ equals the desired probability of drawdown $\phi$, and   the reinsurance-investment strategy is given  in \eqref{eq:strategy1} and \eqref{eq:strategy2}. We complete the proof. \qed
\vspace{3mm}
\begin{remark}\label{re3.2}
We can see from  \eqref{eq:strategy1} and \eqref{eq:strategy2}  that the optimal  strategy is not only depend on the value of surplus $u$ but also on the maximum value  $m$. In particular, similar to Han et al. \cite{HLY19},  the optimal strategy for the drawdown problem coincides with the one for the ruin problem if drawdown does not happen. Setting $r=0$ and $\pi=0$ in our risk model, then   it is easy to check that  the optimal results for the fix ruin level $0$ are identical to those in Bai et al.  \cite{BCZ13}, % (they   studied the ruin probability of a common shock dependence model without investment under the expected-value principle) 
in which  the optimal  retention policy  is independent of the surplus value $u$. 
\end{remark}
When controlling a diffusion towards a goal,  Pestien and Sudderth {\rm \cite{PS85}} showed that,  the optimal strategy which  minimizes  the probability of ruin is the one that maximizes the drift divided by the square of the volatility.  
Moreover, Han et al. {\rm \cite{HLY19}} verified that  the optimal strategy that maximizes the drift divided by the square of the volatility also minimizes the probability of drawdown over an infinite-time horizon.  Thus, it is to be expected from  \eqref{eq:strategy1} and \eqref{eq:strategy2} that  when  the reinsurance premium is calculated by the  expected value principle,  the optimal reinsurance strategy  that maximizes the drift divided by   the square of the volatility is in the form of  excess-of-loss reinsurance, which can be seen  clearly in the following proposition.

\begin{proposition}
The drift of $\hat{U}$ in \eqref{eq:Uhat} divided by the square of its volatility is given by \begin{equation}\label{eq:ratio1}f_1(u,\mH_1,\mH_2, \pi)=\frac{ r u +(\mu-r)\pi+c+\eta(\beta_1+\beta_2)-(1+\eta)(a_1+a_2)}{\sigma^2\pi^2+\gamma_1^2+\gamma^2_2+2 \la\mu_1\mu_2}.\end{equation}
Let ${\mA}_e$ denote  the set of all admissible strategy $\nu^e=(\mathcal{H}_1^e,\mathcal{H}_2^e,\pi)$, where $\mathcal{H}_i^e$ is  excess-of-loss retention  policy   with the following form $$\mathcal{H}^e_{i}(x)=min\{x,d_i\},  ~~~~i=1,2,$$ in which  $0 \leq d_i\leq \infty$. Then, the maximum value of \eqref{eq:ratio1} is attained in set $\mA_e$.  %Obviously,  set $\mA$ contains ${\mA}_e$. 
\end{proposition}
\noindent{\bf Proof:} The proof follows readily from the corresponding proofs given in Borch {\rm \cite{B60}}, Hipp and Taksar {\rm \cite{HT10}}, and Bai et al. \cite{BCZ13}, we omit the details here. \qed
% Borch {\rm \cite{B60}}, Hipp and Taksar {\rm \cite{HT10}} showed that, for any one-dimensional reinsurance policy $\mG$ with $\mG(x)\leq x$ for $x\geq0$ and a nonnegative random variable $Z$, there exists an excess of loss reinsurance $\mG^e(x)=\min\{x,d\}$ with a retention level $0 \leq d\leq \infty$ such that $$\E\left(\mathcal{G}^e(Z)\right)=\int_0^d\pr\{Z>x\}dx=\E\left(\mathcal{G}(Z)\right),$$ $$\E\left(\mathcal{G}^e(Z)\right)^2=\int_0^d2x\pr\{Z>x\}dx\leq\E\left(\mathcal{G}(Z)\right)^2.$$ Hence, let  $\mu^e_i$, $\beta^e_i$  and $\gamma^e_i$ denote the values of $\mu_i$, $\beta_i$ and $\gamma_i$  when $\mathcal{H}_i$  is replaced by $\mathcal{H}_i^e$,  it is easy to find that, for any investment strategy $\pi$, there exists a pair of  excess-of-loss policy $(\mathcal{H}^e_{1}, \mathcal{H}^e_{2})$ such that $$\beta_1+\beta_2=\beta^e_1+\beta^e_2,$$  $$\gamma_1^2+\gamma^2_2+2 \lambda\mu_1\mu_2\geq(\gamma^e_1)^2+(\gamma^e_2)^2+2 \la\mu^e_1\mu^e_2.$$Therefore, we complete the proof.
\vspace{0.5cm}

In the following corollary of Theorem \ref{th3.1}, we consider a special case of our optimization problem by setting $\eta_1=\eta_2$, i.e,  the reinsurance safety loadings for the two classes of business  are the same. We denote the  common reinsurance safety loading by the notation $\eta$.
\begin{corollary}\label{cor3.1} Suppose that $\eta=\eta_1=\eta_2$. Let $g_1$  and $f_1$  be given in Appendix   A.1. The minimum probability of drawdown for the surplus process \eqref{eq:Uhat}  for any $u\in[\alp m,\min(m, u_s)]$ is given by \begin{equation}\label{eq:drawdown4}
\phi(u,m) = \begin{cases}
1 - \dfrac{g_1(u, m)}{g_1(u_s, m)}, &\hbox{if } \alp m \le u \le u_s \le m, \\[5mm]
1 - \displaystyle\exp \left\{ - \int_m^{u_s} f_1(y) dy \right\} \cdot \dfrac{g_1(u, m)}{g_1(u_s, u_s)}, &\hbox{if } \alp m \le u \le m < u_s.
\end{cases}
\end{equation}
The corresponding optimal reinsurance-investment strategy is given in \eqref{eq:strategy1}.  Note that $d^*_1(u)$  is  the  unique solution to  \eqref{eq:G} with  $\eta_i$ $(i=1,2)$   replaced  by $\eta$.
\end{corollary}
\noindent
{\bf Proof:} When $\eta_1=\eta_2$,  $l_X(x)$ and $l_Y(x)$ in \eqref{eq:lxly} are both strictly  increasing functions  with respect to $x\geq0$,   and thus we have  $a_l=0$.  Since $a_k(u)$ is nonnegative, we deduce that  $a_k(u)\geq a_l$ for all  $u\in[\alp m,\min(m, u_s)]$. Following the same lines in Theorem \ref{th3.1}, we obtain the optimal drawdown probability  in \eqref{eq:drawdown4} associated with the optimal strategy in \eqref{eq:strategy1}. We complete the proof.\qed

\subsection{The variance premium principle}
In this subsection, we discuss the  same  optimization problem  under the variance premium principle based on which the reinsurance premium rate can be expressed as 
\begin{equation}\label{eq:var-pre}\delta(\mH_t)=(a_1+a_2-\beta_1-\beta_2)+\Lambda_1\bar{\gamma}_1^2+\Lambda_2\bar{\gamma}^2_2,\end{equation}
where $\Lambda_1$ and $\Lambda_2$ are reinsurer's safety loadings of the two classes of the insurance business,  and $\bar{\gamma}_i$ $(i=1,2)$ is defined as  
\begin{equation}\label{eq:gamma_i}\bar{\gamma}_1=\sqrt{(\la+\la_1)\E\big(X_i-\mathcal{H}_{1}(X_i)\big)^2},~~~~~~\bar{\gamma}_2=\sqrt{(\la+\la_2)\E\big(Y_i-\mathcal{H}_{2}(Y_i)\big)^2}.\end{equation}
Similarly, we plug $\delta(\mH_t)$ above  into \eqref{eq:b_pro} and  define a related  function 
\begin{equation}\label{eq:g2}\begin{array}{lll}g_2(u,\mH_1,\mH_2,\pi) =\big(ru+c+(\mu-r)\pi-a_1-a_2-\Lambda_1\bar{\gamma}_1^2-\Lambda_2\bar{\gamma}^2_2)\xi_u\\[3mm]~~~~~~~~~~~~~~~~~~~~~~~+\displaystyle\frac{1}{2}(\sigma^2\pi^2+\gamma^2_1+\gamma^2_2+2\la\mu_1\mu_2)\xi_{uu}.\end{array}\end{equation}
By using the cumulative distribution functions of $X$ and $Y$, respectively,  we  rewrite $g_2$ as follows:
\begin{equation}\label{eq:ng2}\begin{array}{lll}g_2(u,\mathcal{H}_1,\mathcal{H}_2,\pi) =\big( ru-\kappa+(\mu-r)\pi
-\Lambda_2(\la+\la_2)(\E\mH^2_2-2\E(\mH_2 Y))\big)\xi_u
\\[3mm]~~~~~~~~~~~~~~~~~~~~~~~+\dfrac{1}{2}\big(\sigma^2\pi^2+(\la+\la_2)\E\mH^2_2\big)\xi_{uu}+\displaystyle\int_0^\infty\Big\{-\Lambda_1(\la+\la_1)\big(\mH_1^2(x)-2x\mH_1(x)\big)\xi_u\\[5mm]~~~~~~~~~~~~~~~~~~~~~~~+\displaystyle\frac{1}{2}\big((\la+\la_1)\mH^2_1(x)+2\la\mH_1(x)\E\mH_2\big)\xi_{uu}\Big\}dF_X(x)\\[3mm]~~~~~~~~~~~~~~~~~~~~=\big( ru-\kappa+(\mu-r)\pi-\Lambda_1(\la+\la_1)(\E\mH^2_1-2\E(\mH_1 X))\big)\xi_u\\[3mm]~~~~~~~~~~~~~~~~~~~~~~~ +\dfrac{1}{2}\big(\sigma^2\pi^2+(\la+\la_1)\E\mH^2_1\big)\xi_{uu}+\displaystyle\int_0^\infty\Big\{-\Lambda_2(\la+\la_2)(\mH_2^2(x)-2x\mH_2(x))\xi_u\\[3mm]~~~~~~~~~~~~~~~~~~~~~~~+\displaystyle\frac{1}{2}\big((\la+\la_2)\mH^2_2(x)+2\la\mH_2(x)\E\mH_1\big)\xi_{uu}\Big\}dF_Y(x).\end{array}
\end{equation}
%which we can minimize $x$-by-$x$ subject to the condition  $0 \le \mH_i(x) \le x$ $(i=1,2)$. 
Therefore, the minimizer of $g_2$ is obtained  at 
\begin{equation}\label{eq:cand_H*}\mH^*_1(u, x)=q^*_1(u,x) \wedge x,~~~~~\mH^*_2(u, y)= q^*_2(u,x)\wedge y\end{equation}
with
$$
q^*_1(u,x)=\displaystyle\frac{2\Lambda_1}{2\Lambda_1-\frac{\xi_{uu}}{\xi_u}}x+\frac{\la\E\mH_2\frac{\xi_{uu}}{\xi_u}}{(\la+\la_1)(2\Lambda_1-\frac{\xi_{uu}}{\xi_u})},$$
$$
q^*_2(u,x)=\displaystyle\frac{2\Lambda_2}{2\Lambda_2-\frac{\xi_{uu}}{\xi_u}}y+\frac{\la\E\mH_1\frac{\xi_{uu}}{\xi_u}}{(\la+\la_2)(2\Lambda_2-\frac{\xi_{uu}}{\xi_u})},$$
and \begin{equation}\label{eq:cand_pi*}\pi^*(u)=-\displaystyle\frac{u-r}{\sigma^2}\frac{\xi_u}{\xi_{uu}}.\end{equation}
%Here, even though  we can find a candidate extreme minimum point  of  $g_2$, 
With some tedious calculations, 
we find that  the uniqueness and existence of   the extreme minimum point and the  corresponding  value function are not easy to determined when the reinsurance premium is calculated by \eqref{eq:var-pre}.
However, inspired by Han et al. \cite{HLY18}, if the reinsurance safety loadings for the  two classes are   the same, it may be possible to derive explicit expressions for the optimal results. Therefore, in the rest of this subsection, we focus on investigating the optimization problem with a common reinsurance safety loading.

When the reinsurance premium is calculated by the variance premium principle with a common safety loading $\Lambda$, we have 
\begin{equation}\label{eq:var-pre2}\delta(\mH_t)=(a_1+a_2-\beta_1-\beta_2)+\Lambda(\bar{\gamma}_1^2+\bar{\gamma}^2_2+2\la\bar{\mu}_1\bar{\mu}_2),\end{equation}
where $\bar{\gamma}_i$ $(i=1,2)$ is given in \eqref{eq:gamma_i}, and $\bar{\mu}_i$ $(i=1,2)$ is defined as  
\begin{equation}\label{eq:mu_i}\bar{\mu}_1=\E\big(X_i-\mathcal{H}_{1}(X_i)\big),~~~~~~\bar{\mu}_2=\E\big(Y_i-\mathcal{H}_{2}(Y_i)\big).\end{equation}
We define the function $\bar{g}_2$ as follows \begin{equation}\label{eq:barg2}\begin{array}{lll}\bar{g}_2(u,\mH_1,\mH_2,\pi) =\big(ru+c+(\mu-r)\pi-a_1-a_2-\Lambda(\bar{\gamma}_1^2+\bar{\gamma}^2_2+2\la\bar{\mu}_1\bar{\mu}_2)\big)\xi_u\\[5mm]~~~~~~~~~~~~~~~~~~~~~~~~+\displaystyle\frac{1}{2}(\sigma^2\pi^2+\gamma^2_1+\gamma^2_2+2\la\mu_1\mu_2)\xi_{uu}.\end{array}\end{equation}
%\begin{equation}\label{eq:ng2}\begin{array}{lll}g_2(u,\mathcal{H}_1,\mathcal{H}_2,\pi) =\big( ru+c-\delta(0)+(\mu-r)\pi-\Lambda\big((\la+\la_2)(\E\mH^2_2-2\E(\mH_2 Y))-2\la\E X\E \mH_2\big)\big)h_u\\[5mm]~~~~~~~~~~~~~~~~~~~~~~~ +\dfrac{1}{2}\big(\sigma^2\pi^2+(\la+\la_2)\E\mH^2_2\big)h_{uu}+\displaystyle\int_0^\infty\Big\{-\Lambda\big((\la+\la_1)\mH_1^2(x)-2(\la+\la_1)x\mH_1(x)\\[5mm]~~~~~~~~~~~~~~~~~~~~~~~+2\la\mH_1(x)\E\mH_2-2\la\mH_1(x)\E Y\big)h_u+\displaystyle\frac{1}{2}\big((\la+\la_1)\mH^2_1(x)+2\la\mH_1(x)\E\mH_2\big)h_{uu}\Big\}dF_X(x).\\[5mm]~~~~~~~~~~~~~~~~~~~~=\big( ru+c-\delta(0)+(\mu-r)\pi-\Lambda\big((\la+\la_1)(\E\mH^2_1-2\E(\mH_1 X))-2\la\E Y\E \mH_1\big)\big)h_u\\[5mm]~~~~~~~~~~~~~~~~~~~~~~~ +\dfrac{1}{2}\big(\sigma^2\pi^2+(\la+\la_1)\E\mH^2_1\big)h_{uu}+\displaystyle\int_0^\infty\Big\{-\Lambda\big((\la+\la_2)\mH_2^2(x)-2(\la+\la_2)x\mH_2(x)\\[5mm]~~~~~~~~~~~~~~~~~~~~~~~+2\la\mH_2(x)\E\mH_1-2\la\mH_2(x)\E X\big)h_u+\displaystyle\frac{1}{2}\big((\la+\la_2)\mH^2_2(x)+2\la\mH_2(x)\E\mH_1\big)h_{uu}\Big\}dF_Y(x)\end{array}\end{equation}
%From this integral representation of $g_2$, we deduce that we can minimize $g_2$ by mimizing the integrand $x$-by-$x$, subject to $0 \le \mH_1(x) \le x$.  As a function of $\mH_1(x)$, the integrand is a parabola, so it is minimized  by
It turns out that  the minimizer of $\bar{g}_2$ is obtained  at 
$$\mH^*_1(u,x)=\bar  q^*_1(u,x) \wedge x,~~~~~~~\mH^*_2(u,y)= \bar  q^*_2(u,y)\wedge y$$  
with
\begin{equation}\label{eq:q1}
\bar q^*_1(u,x)=\displaystyle-\frac{\la\E\mH_2}{\la+\la_1}+\frac{2\Lambda}{2\Lambda-\frac{\xi_{uu}}{\xi_u}}x+\frac{2\Lambda\la\E Y}{(\la+\la_1)(2\Lambda-\frac{\xi_{uu}}{\xi_u})},\end{equation}
\begin{equation}
\label{eq:q2}\bar  q^*_2(u,y)=\displaystyle-\frac{\la\E\mH_1}{\la+\la_2}+\frac{2\Lambda}{2\Lambda-\frac{\xi_{uu}}{\xi_u}}y+\frac{2\Lambda\la\E X}{(\la+\la_2)(2\Lambda-\frac{\xi_{uu}}{\xi_u})},\end{equation}
and $\pi^*(u)$ is given in the form of \eqref{eq:cand_pi*}. %$$\pi^*(u)=-\displaystyle\frac{u-r}{\sigma^2}\frac{\xi_u}{\xi_{uu}}.$$
%{\bf \color{red} Obviously, the forms of $\mH_i(x)$ are not the proportional reinsurance. Once we  obtain the ratio of $\frac{h_{uu}}{h_{u}}$,  the drawdown probability can be calculated!} 
Firstly, we suppose that  $(\mH^*_1, \mH^*_2)=(q_1^*,q_2^*)$, and then taking the expected value on both sides of  equations  \eqref{eq:q1} and \eqref{eq:q2} yields 
$$\E \mH_1=\frac{2\Lambda}{2\Lambda-\frac{\xi_{uu}}{\xi_u}}\E X+\frac{2\Lambda\la}{(\la+\la_1)(2\Lambda-\frac{\xi_{uu}}{\xi_u})}\E Y-\frac{\la}{\la+\la_1}\E\mH_2,$$
$$\E \mH_2=\frac{2\Lambda}{2\Lambda-\frac{\xi_{uu}}{\xi_u}}\E Y+\frac{2\Lambda\la}{(\la+\la_2)(2\Lambda-\frac{\xi_{uu}}{\xi_u})}\E X-\frac{\la}{\la+\la_2}\E\mH_1.$$
By solving the equations above, it gives 
\begin{equation}\label{eq:H}\E\mH_1=\frac{2\Lambda}{(2\Lambda-\frac{\xi_{uu}}{\xi_u})} \E X,~~~~~~\E\mH_2=\frac{2\Lambda}{(2\Lambda-\frac{\xi_{uu}}{\xi_u})} \E Y.\end{equation}
Bringing  \eqref{eq:H} back  into \eqref{eq:q1}  and \eqref{eq:q2}, it then follows that \begin{equation}\label{eq:qstar}\bar  q^*_i=\displaystyle\frac{2\Lambda}{2\Lambda-\frac{\xi_{uu}}{\xi_u}}x, ~~~~~~~~i=1,2.\end{equation}
For $\frac{\xi_{uu}}{\xi_{u}}<0$, it is easy to check that $\bar q^*_i(u,x)<x$, and hence  we indeed  have
 $\mH^*_i=\bar q^*_i~(i=1,2).$ Besides,  the Hessen matrix of  $\bar{g}_2$  at point $(\mH^*_1,\mH^*_2,\pi^*)$ can be written as  
\begin{equation}\label{eq:matrixDbar}\bar{\textbf{D}}=\left(\begin{array}{ccccc}
&(\la+\la_1)\E X^2(\xi_{uu}-2\Lambda \xi_u)& \la \E X\E Y(\xi_{uu}-2\Lambda \xi_u) &0\\[3mm]
& \la \E X\E Y(\xi_{uu}-2\Lambda \xi_u)&(\la+\la_2)\E Y^2(\xi_{uu}-2\Lambda \xi_u) &0\\[3mm]
&0&0 & \sigma^2\xi_{uu}
 \end{array}\right),\end{equation}
which is   positive definite. Therefore, if we can find a convex candidate solution $\xi$, then the point $(\mH^*_1,\mH^*_2,\pi^*)=(q^*_1,q^*_2,\pi^*)$ is the extreme minimum point of $\bar{g}_2$.
%$$\frac{\xi_{uu}}{\xi_u}=\frac{2\Lambda (ru-\kappa)+A+2\Lambda^2C+\sqrt{(2\Lambda (ru-\kappa)+A+2\Lambda^2C)^2-8\Lambda A(ru-\kappa)}}{2(ru-\kappa)}<0$$
Before we continue to track our value function, we first give a lemma which plays a similar role as Lemma \ref{lem:4.3}.

For notational convenience, we denote  \begin{equation}\label{eq:A1} A_1=(\la+\la_1)\E X^2+(\la+\la_2)\E Y^2+2\la\E X\E Y.\end{equation} 
\begin{lemma}\label{lem:4.4} Recall   $\mM_t$  of  \eqref{eq:mM_t} and $b \in (\alpha m, u_s)$. Let $\delta(\mH_t)$  and $A_1$ be given in \eqref{eq:var-pre2} and  \eqref{eq:A1}, respectively.   We define  that     $\Delta=\frac{(\mu-r)^2}{2\sigma^2}$. Then, there exists $\vartheta=\tilde{\gamma}$ such that $\mM_t\geq\frac{\kap-rb}{2}$ for all $0\leq t\leq \tau_{\alpha b}$ with  $\tilde{\gamma}$  given by \begin{equation}\label{eq:wgam}{\tilde{\gamma}}=-\Lambda-\displaystyle\frac{\Delta+2\Lambda^2A_1+\sqrt{(\Lambda (rb-\kappa)+\Delta+2\Lambda^2A_1)^2-4\Lambda \Delta(rb-\kappa)}}{rb-\kappa}\in(0,\infty)\end{equation}\end{lemma}
\begin{proof}   Instituting  \eqref{eq:var-pre2}  back into the expression of $\mM_t$ yields
 $$\mM_t=-r\hU_t-c-(\mu-r)\pi_t+a_1+a_2+\Lambda(\bar{\gamma}_1^2+\bar{\gamma}^2_2+2\la\bar{\mu}_1\bar{\mu}_2)+\frac{1}{2}\vartheta(\sigma^2\pi_t^2+\gamma^2_1+\gamma^2_2+2\la\mu_1\mu_2).$$
%We wish to show that, for $\delta$ large enough, $\mM_t\geq\frac{\kappa-rb}{2}>0$ for all $t\geq0$ with probability 1. 
Since the inequalities $\kap-r\hU_t\geq \kap-rb>0$ hold for all $0\leq t\leq \tau_{\alpha b}$ , we  may  find $\vartheta>0$ such that 
$$\begin{array}{ll}%\Lambda\left((\la+\la_1)\E(\mH_1^2-2\mH_1X)+(\la+\la_2)\E(\mH_2^2-2\mH_2Y)+2\la\E(\mH_1\mH_2-\mH_1Y-\mH_2X)\right)
\Lambda(\bar{\gamma}_1^2+\bar{\gamma}^2_2+2\la\bar{\mu}_1\bar{\mu}_2)-\delta(\textbf{0})-(\mu-r)\pi+\displaystyle\frac{1}{2}\vartheta(\sigma^2\pi^2+\gamma^2_1+\gamma^2_2+2\la\mu_1\mu_2)\geq-\frac{\kap-rb}{2},\end{array}$$
or equivalently, \begin{equation}\label{eq:var_ine1}\frac{1}{2}\vartheta\geq\frac{%-\Lambda\left((\la+\la_1)\E(\mH_1^2-2\mH_1X)+(\la+\la_2)\E(\mH_2^2-2\mH_2Y)+2\la\E(\mH_1\mH_2-\mH_1Y-\mH_2X)\right)
-\Lambda(\bar{\gamma}_1^2+\bar{\gamma}^2_2+2\la\bar{\mu}_1\bar{\mu}_2)+\delta(\textbf{0})+(\mu-r)\pi-\frac{\kap-rb}{2}}{\sigma^2\pi^2+\gamma^2_1+\gamma^2_2+2\la\mu_1\mu_2}\end{equation}
for all admissible strategies $\nu=(\mH_1,\mH_2,\pi)$.
In other words, we wish to show that the right side of  inequality \eqref{eq:var_ine1} has a finite maximum over such $\nu$.  By comparing the right side of the inequality with $\bar{g}_2$ in \eqref{eq:barg2}, we deduce that the right side of the inequality is maximized by
$$\tilde{\mH}_1(x)=\displaystyle-\frac{\la\E\mH_2}{\la+\la_1}+\frac{2\Lambda}{2\Lambda+\tilde{\gamma}}x+\frac{2\Lambda\la\E Y}{(\la+\la_1)(2\Lambda+\widetilde{\gamma})},$$
$$\widetilde{\mH}_2(y)=\displaystyle-\frac{\la\E\mH_1}{\la+\la_2}+\frac{2\Lambda}{2\Lambda+\widetilde{\gamma}}y+\frac{2\Lambda\la\E X}{(\la+\la_2)(2\Lambda+\widetilde{\gamma})},$$  
and $$\tilde{\pi}=\displaystyle\frac{1}{\tilde{\gamma}}\cdot\frac{u-r}{\sigma^2},$$
in which $\tilde{\gamma}$ uniquely solves 
$$-\Lambda(\bar{\gamma}_1^2+\bar{\gamma}^2_2+2\la\bar{\mu}_1\bar{\mu}_2)+\delta(\textbf{0})+(\mu-r)\pi-\tilde{\gamma}(\sigma^2\pi^2+\gamma^2_1+\gamma^2_2+2\la\mu_1\mu_2)=\frac{\kap-rb}{2}.$$ 
With some calculation and rearrangement,  we obtain $\tilde{\gamma}$ given in \eqref{eq:wgam}. By substituting $(\tilde{\mH}_1,\tilde{\mH}_2,\tilde{\pi})$  into the right side of \eqref{eq:var_ine1}, we obtain  $\frac{\vartheta}{2}>\frac{\tilde{\gamma}}{2}$; thus, set $\vartheta=\tilde{\gamma}\in(0,\infty),$ and for this choice of $\vartheta$, we have $\mM_t\geq\frac{\kap-rb}{2}>0$ for all $0\leq t\leq\tau_{\alpha b}$ with probability 1. We complete the proof.
\end{proof}
In  the following theorem,  combining with the verification theorem given in Section \ref{sec:3}, we present the solution to our optimization problem under the variance premium principle with a common safety loading.
\begin{theorem}\label{thm:4.2}
 Let $\bar{g}_1(u,m)$ and $\bar{f}_1(y)$ be given in Appendix  A.2. Let  $q^*_1$, $q^*_2$ and $\pi^*$  be given in  \eqref{eq:cand_pi*}, \eqref{eq:q1} and \eqref{eq:q2}, respectively. 
 The minimum probability of drawdown for the surplus process \eqref{eq:Uhat} on $\mO$ is given by
\begin{equation}\label{eq:drawdown5}
 \phi(u,m)=\left\{\begin{array}{llll}
 \displaystyle1-\frac{\bar{g}_1(u,m)}{\bar{g}_1(u_s,m)},& \alpha m\leq u \leq u_s\leq m, \\[5mm]
 \displaystyle1- \exp\left\{-\int_m^{u_s}{\bar{f}_1(y)dy}\right\}\cdot\frac{\bar{g}_1(u,m)}{\bar{g}_1(u_s,u_s)}, & \alpha m\leq u \leq m < u_s.
 \end{array}\right.
\end{equation}
Also,  the corresponding  optimal reinsurance-investment strategy, for all $t\geq0$, is given in feedback form by
 \begin{equation}\label{eq:strategy3}\left \{\begin{array}{ll}\mH^*_{1t}=\displaystyle\frac{2\Lambda}{2\Lambda-{\chi}_1(\hat{U}_t)}X,~~~~~~
\mH^*_{2t}=\displaystyle\frac{2\Lambda}{2\Lambda-{\chi}_1(\hat{U}_t)}Y,\\[5mm]
\pi^*_t=\displaystyle\frac{\mu-r}{\sigma^2}\cdot\frac{1}{{\chi}_1(\hat{U}_t)},
\end{array}\right.\end{equation} in which  $\chi_1(u)$ is given by \begin{equation}\label{eq:beta1}\chi_1(u)=\Lambda +\displaystyle\frac{\Delta+2\Lambda^2A_1+\sqrt{(2\Lambda (ru-\kappa)+\Delta+2\Lambda^2A_1)^2-8\Lambda \Delta(ru-\kappa)}}{2(ru-\kappa)}.\end{equation}
  \end{theorem}
  \noindent{\bf \emph{Proof}:}  Substituting  $(\mH^*_1,\mH^*_2,\pi^*)=(\bar q^*_1,\bar q^*_2,\pi^*)$ back into  \eqref{eq:b_pro}, one can show that 
\begin{equation}\label{eq:ratio-var}ru-\kappa-\displaystyle \Delta\cdot\frac{\xi_u}{\xi_{uu}}+A_1\cdot\frac{2\Lambda^2}{2\Lambda-\frac{\xi_{uu}}{\xi_u}}=0\end{equation}with 
%\begin{equation}\label{eq:A1} A_1=(\la+\la_1)\E X^2+(\la+\la_2)\E Y^2+2\la\E X\E Y.\end{equation}
$A_1$ be given in \eqref{eq:A1}. Solving the ratio $\frac{\xi_{uu}}{\xi_u}$   from the equation \eqref{eq:ratio-var} yields 
   $\frac{\xi_{uu}}{\xi_u}=\chi_1(u)<0$.
Let $\xi$ equal to the right-hand side of \eqref{eq:drawdown5}, following the similar arguments and steps in Theorem \ref{th3.1}, we can verify that $\xi$  satisfies all the conditions stated in Theorem \ref{thm:verif}. Finally, we have $\phi=\xi$ with the optimal reinsurance-investment  strategy  given in \eqref{eq:strategy3}. $\hfill\Box$
\begin{remark}\label{rem4.1} %The drift of $\hat{U}^\nu$ divided by the square of its volatility is given by \begin{equation}\label{eq:ratio2}f_2(u,\mH_1,\mH_2,\pi)=\frac{ r u +(\mu-r)\pi+c-a_1-a_2-\Lambda(\bar{\gamma}_1^2+\bar{\gamma}^2_2+2\la\bar{\mu}_1\bar{\mu}_2)}{\sigma^2\pi^2+\gamma_1^2+\gamma^2_2+2 \la\mu_1\mu_2}.\end{equation} Let ${\mA}_p$ denote  the set of all admissible strategy $\nu^p=(\mathcal{H}_1^p,\mathcal{H}_2^p,\pi)$, where $\mathcal{H}_i^p$ is  proportional reinsurance policy   with the following form $$\mathcal{H}^p_{i}(x)=q_i\,x,  ~~~~i=1,2,$$ where $0 \leq q_i\leq 1$. % Obviously,  set $\mA$ also contains ${\mA}_p$. Then, we have the following proposition. For a variance premium  principle and  any admissible policy $(\mathcal{H}_1,\mathcal{H}_2,\pi)\in{\mA}$,  the maximum value in \eqref{eq:ratio2} is attained in set ${\mA}_p$.
 Hipp and Taksar {\rm \cite{HT10}} showed that,  for any one-dimensional reinsurance policy $\mG$ with $\mG(x)\leq x$ for $x\geq0$, %and a nonnegative random variable $Z$, 
  there exists a proportional reinsurance  policy $\mG^p(x)=q\,x$ with a risk sharing proportional  $0 \leq q\leq 1$ such that  the ratio of the drift divided the square of the volatility is maximized  under the variance premium principle.
%$$\E\left(\mathcal{G}^e(Z)\right)=\int_0^d\pr\{Z>x\}dx=\E\left(\mathcal{G}(Z)\right),$$ $$\E\left(\mathcal{G}^p(Z)\right)^2=q^2\,\E\big(Z^2\big)=\E\left(\mathcal{G}(Z)\right)^2,$$$$\E\left(Z-\mathcal{G}^p(Z)\right)^2=\E\big(Z^2\big)-2q\E\big(Z^2\big)+q^2\,\E\big(Z^2\big)\leq\E\left(Z-\mathcal{G}(Z)\right)^2,$$
%Hence, let  $\mu^p_i$, $\bar{\mu}^p_i$, $\gamma^p_i$  and $\bar{\gamma}^p_i$ denote the values of $\mu_i$, $\bar{\mu}_i$, $\gamma_i$ and $\sigma_i$   when $\mathcal{H}_i$  is replaced by $\mathcal{H}_i^p$,  it is easy to find that there exists a pair of  proportional reinsurance policy $(\mathcal{H}^p_{1}, \mathcal{H}^p_{2})$ such that 
%$$(\bar{\gamma}^p_1)^2+(\bar{\gamma}^p_2)^2<\gamma^2_1+\gamma^2_2,$$  $$(\gamma^p_1)^2+(\gamma^p_2)^2=\gamma_1^2+\gamma^2_2.$$
When we consider  the risk model with dependence structure, it is not easy to  obtain the optimal reinsurance strategy by the same method in Hipp and Taksar {\rm \cite{HT10}} because of  the  item $2\la\mu_1\mu_2 $ in volatility, but  the optimal reinsurance strategy is still in the form of quota-share reinsurance, which we can see  clearly from \eqref{eq:strategy3}.   \end{remark}
\section{Further discussions}\label{sec:5}
In this section, we may extend our work to the  case where the  insurance company involves $n$ $(n\geq3)$ dependent classes of insurance business.  Inspired by the analysis in Section \ref{sec:4},  the optimal results can also be derived  explicitly  not only for the expected value principle (under  some  hypothesis of the parameters) but also for the variance premium principle.
  
  Let $\{X^{(l)}_i, i\geq 1\}$ $(l=1,2,...,n)$ be the claim size random variables for the $l$th  class of business following a common distribution $F_{X^{(l)}}(x)$. Assume that  $F_{X^{(l)}}(x)=0$ for $x\leq 0$,  $0<F_{X^{(l)}}(x)\leq1$ for $x>0$. Then,  
the aggregate claims up to time $t$ for the $l$th line of business are denoted by
$\sum_{i=1}^{N_l(t)+N(t)}X^{(l)}_i$, where $N_l(t)+N(t)$  is the claim number process for class $l$
 $(l=1,2,...,n)$. Assume that    $  N_l(t)$ and $N(t)$  are $n+1$ independent Poisson processes with parameters $\lambda_1,... ,\lambda_n$ and  $\lambda$, respectively. 
 %\subsection{ Some results for the expected-value principle}
 Let $\mH=(\mathcal{H}_1,\mathcal{H}_{2},..., \mH_n)$ and $\pi$ denote the  retention level and the investment strategy, respectively.  For convenience, we assume that the reinsurance safety loadings for all classes are the same in the following context.  Under the expected value principle,  the reinsurance premium rate  at time $t$ is $$\delta(\mH_t)=(1+\eta)\sum_{l=1}^{n}(a_l-\beta_l)$$ 
 with 
$$ a_l=(\la+\la_l)\E X^{(l)}_i,~~~~ \beta_l=(\la+\la_l)\E\big(\mathcal{H}_{l}(X^{(l)}_i)\big).$$
We point out that when  the insurer involves  $n$ dependent classes of insurance business, it is not easy to obtain the optimal results following the same lines as in Section \ref{sec:4.1}.  In fact,  we need to define $n$ numbers of   auxiliary functions similar to \eqref{eq:lxly}, but the monotonicity of the functions can not  be determined clearly.   Here, to keep things simple, we constrain the retention strategy in the interval $[0,\infty)$. On the one hand, for $\mH_l(x) \in[0,x]~(l=1,2,\cdots,n)$, the insurer has a share of
 reinsurance cover. On the other hand, the case with $\mH_l(x) \in(x, \infty)$ may be thought of as acquiring
new business. We still assume   $\xi$ is  the candidate solution to optimization problem. Then  the differential operator $\mA^\nu$ in \eqref{eq:op} can be rewritten as
\begin{equation}\label{eq:op1}
\mA^{\nu} \xi(u,m) = \Big[ r u -\kappa+(\mu-r)\pi_t+\eta\sum_{l=1}^n\beta_l\Big] \xi_u + \frac{1}{2}\Big(\sigma^2\pi^2+\sum_{l=1}^n\gamma_l^2+\sum _{k\neq j} \la\mu_k\mu_j\Big) \xi_{uu},
 \end{equation}
where 
$$\mu_l=\E\mathcal{H}_{l}\big(X^{(l)}_i\big),~~~~\gamma_l=\sqrt{(\la+\la_l)\E\left(\mathcal{H}_{l}(X^{(l)}_i)\right)^2}$$
 for  $l, j, k=1,2,\ldots, n.$   By differentiation, the minimizer of $\mA^{\nu} \xi(u,m)$  can be  obtained at
 \begin{equation}\label{eq:Hn}\mH^*_l (u)=-\eta\frac{\xi_u}{\xi_{uu}}-\frac{\la}{\la+\la_l}\sum_{j\neq l}\E\mH^*_j(X^{(j)}_i)\end{equation} for $l=1,2,\cdots,n,$ and  \begin{equation}\label{eq:pin}\pi^*(u)=-\displaystyle\frac{u-r}{\sigma^2}\frac{\xi_u}{\xi_{uu}}.\end{equation} Evaluating the expected value on both sides  of the equations \eqref{eq:Hn}  for each $l$ yields \begin{equation}\label{eq:EHi}\textbf{A}(\E\mH_1,\E\mH_2,\cdots,\E\mH_n)^T=-\eta\frac{\xi_u}{\xi_{uu}}\textbf{1}^T,\end{equation}
 where the matrix  $\textbf{A}$ is given as
  \begin{equation}\label{eq:matrixA}\textbf{A}=\left(\begin{array}{ccccc}
 &1&\frac{\la}{\la+\la_1}& \cdots & \frac{\la}{\la+\la_1}\\[3mm]
 &\frac{\la}{\la+\la_2}&1& \cdots& \frac{\la}{\la+\la_2}\\[3mm]
 &\vdots&\vdots&\ddots&\vdots\\[3mm]
 &\frac{\la}{\la+\la_n}&\frac{\la}{\la+\la_n} &\cdots& 1
  \end{array}\right),\end{equation}
  and  $\textbf{1}=(1,1,\ldots,1).$
Besides,   for any $x=(x_1,x_2,\cdots,x_n)\in\R^n$ and $x\neq 0$, we have
$$ \begin{array}{ll}x\cdot\textbf{A}\cdot x^T&=\displaystyle\frac{1}{(\la+\la_1)\cdots(\la+\la_n)}\Big(\sum_{i=1}^n(\la+\la_i)x^2_i+\la\sum_{i\neq j}x_ix_j\Big)\\[5mm]&=\displaystyle\frac{\la(x_1+x_2+\cdots+x_n)^2}{(\la+\la_1)\cdots(\la+\la_n)}+\sum_{i=1}^n\la_ix^2_i>0,\end{array}$$
 which implies that $\textbf{A}$ is a positive definite matrix.  Therefore, it follows that 
 $$(\E\mH_1,\E\mH_2,\cdots,\E\mH_n)^T=-\eta\frac{\xi_u}{\xi_{uu}}\textbf{A}^{-1}\textbf{1}^T,$$
from which we can get
\begin{equation}\label{eq:sumEHi}\sum_{l=1}^{n}\E\mH_l=-\eta\frac{\xi_u}{\xi_{uu}}\textbf{1}\textbf{A}^{-1}\textbf{1}^T.\end{equation}
It is clear  from \eqref{eq:Hn} that  $\mH^*_l (u)=\E\mH^*_l (u)$,  then combining with the result in \eqref{eq:sumEHi} yields \begin{equation}\label{eq:Histar}\mH^*_l (u)=-\eta K_l\frac{\xi_u}{\xi_{uu}}\end{equation}
with \begin{equation}\label{eq:K_l}K_l=\frac{\la+\la_l}{\la_l}-\textbf{1}\textbf{A}^{-1}\textbf{1}^T.\end{equation}
Denote that \begin{equation}\label{eq:barDelta}\bar\Delta=\frac{1}{2}\left(\frac{\mu-r}{\sigma}\right)^2+\frac{1}{2}\eta^2\left(\sum_{l=1}^n(\la+\la_l)(2K_l-K^2_l)-\la\sum_{j\neq k}K_jK_k\right).\end{equation}
\begin{assumption}  We assume that $\bar\Delta$ defined in \eqref{eq:barDelta} is positive.\end{assumption}
 For $\frac{\xi_{uu}}{\xi_{u}}<0$, we can see from \eqref{eq:Histar} that $\mH^*_l (u)>0$   if  the inequality
\begin{equation}\label{eq:constraint}\textbf{1}\textbf{A}^{-1}\textbf{1}^T<\frac{\la+\la_l}{\la}\end{equation} holds  for all $j, l=1,2,\cdots, n$. Once the inequality \eqref{eq:constraint} holds, bringing \eqref{eq:pin} and  \eqref{eq:Histar}  back into \eqref{eq:op1} and putting $\mA^{\nu} \xi(u,m)=0$ yield
\begin{equation}\label{eq:ratio3}\frac{\xi_{uu}}{\xi_u}=\frac{\bar\Delta}{ru-\kappa}<0.\end{equation}
\begin{theorem}\label{thm:5.1}
 Suppose $\bar\Delta>0$ and the inequality \eqref{eq:constraint} holds. Let ${g}_3(u,m)$ and ${f}_3(y)$ be given in Appendix  A.1.  Let $K_l~(l=1,2,\dots,n)$ be given in \eqref{eq:K_l}. The minimum probability of drawdown for $n$ dependent classes of insurance business on $\mO$ is given by
$$
 \phi(u,m)=\left\{\begin{array}{llll}
 \displaystyle1-\frac{{g}_3(u,m)}{{g}_3(u_s,m)},& \alpha m\leq u \leq u_s\leq m, \\[5mm]
 \displaystyle1- \exp\left\{-\int_m^{u_s}{{f}_3(y)dy}\right\}\cdot\frac{{g}_3(u,m)}{{g}_3(u_s,u_s)}, & \alpha m\leq u \leq m <u_s.
 \end{array}\right.
$$
Also,  the corresponding  optimal reinsurance-investment strategy, for all $t\geq0$ and  l=1,2,\dots,n, is given in feedback form by
 \begin{equation}\label{eq:strategy4}%\left \{\begin{array}{ll}
 \mH^*_{lt}=-\displaystyle\frac{\eta K_l(r\hat{U}_t-\kappa)}{\bar\Delta},
~~~~\pi^*_t=-\displaystyle\frac{(\mu-r)(r\hat{U}_t-\kappa)}{\sigma^2\bar\Delta}.
\end{equation}
 \end{theorem}
\begin{remark}
 %When we  extend our work to the  case where the  insurance company involves $n$ dependent classes of insurance business under the expected-value principle,  to keep things simple,  we constrain the reinsurance strategy in the interval $[0,\infty)$.    Under the assumption of \eqref{eq:constraint},  we deduce that the optimal reinsurance strategy is  in the form of \eqref{eq:Histar} for $\frac{\xi_u}{\xi_{uu}}<0$.  
 Provided that the condition \eqref{eq:constraint} holds, we conjecture that the ratio given in \eqref{eq:ratio3} should be negative under the strategy in \eqref{eq:Histar} which minimizes the probability of drawdown. In fact,  we have confirmed by some numerical examples that this conjecture is correct, and  invite the interested reader to prove this  hypothesis theoretically. Here, we suppose that $\bar\Delta>0$. \end{remark}

When the reinsurance premium is calculated according to the variance premium principle,  the reinsurance premium rate at time $t$   becomes 
 $$\delta(\mathcal{H}_t)=\sum_{l=1}^n(a_l-\beta_l)+\Lambda\Big(\sum_{l=1}^n\bar{\gamma}_l^2+\la\sum_{j\neq k}\bar{\mu}_j \bar{\mu}_k\Big) $$
 with 
$$\bar{\mu}_l=\E\big(X^{(l)}_i-\mathcal{H}_{l}(X^{(l)}_i)\big),~~~~~\bar{\gamma}_l=\sqrt{(\la+\la_l)\E\big(X^{(l)}_i-\mathcal{H}_{l}(X^{(l)}_i)\big)^2}$$
 for  $j, k,l =1,2,\ldots, n.$  As a result,  the differential operator $\mA^\nu$ in \eqref{eq:op} can be rewritten as
  \begin{equation}\label{eq:op2}\begin{array}{lll}\mA^{\nu} \xi(u,m) =\Big[ru+c+(\mu-r)\pi-\sum_{l=1}^na_l-\Lambda\big(\sum_{l=1}^n\bar{\gamma}_l^2+\la\sum_{j\neq k}\bar{\mu}_j \bar{\mu}_k\big)\Big]\xi_u\\[5mm]~~~~~~~~~~~~~~~~~+\displaystyle\frac{1}{2}\big(\sigma^2\pi^2+\sum_{l=1}^n\gamma_l^2+\sum _{k\neq j} \la\mu_k\mu_j\big)\xi_{uu}.
\end{array}\end{equation}
   To  get the optimal  reinsurance strategy $\mH^*=(\mH^*_1,\mH^*_{2},..., \mH^*_n) $, we need to  solve the following system  \begin{equation}\label{eq:EV}(\E \mH_1, \E \mH_2,  \cdots, \E \mH_n)\cdot\textbf{A}^T=\frac{2\Lambda}{2\Lambda-\frac{\xi_{uu}}{\xi_u}}\textbf{D}^T,\end{equation}
in which  the vector $\textbf{D}$ is given as 
 $$\textbf{D}=\left(\begin{array}{ccccc}
&\E X^{(1)}+\frac{\la}{\la+\la_1}\sum_{l\neq1}\E X_i^{(l)}\\[3mm]&\E X^{(2)}+\frac{\la}{\la+\la_2}\sum_{l\neq2}\E X_i^{(l)}\\[3mm]&\vdots\\[3mm]&\E X^{(n)}+\frac{\la}{\la+\la_n}\sum_{l\neq n}\E X_i^{(l)}
\end{array}\right),$$
and matrix $\textbf{A}$ is given in \eqref{eq:matrixA}. Inspired  by the optimal results for $n=2$, it is not difficult to verify that
 \begin{equation}\label{eq:solEV}\E \mH_l=\frac{2\Lambda}{2\Lambda-\frac{\xi_{uu}}{\xi_u}}\E X_i^{(l)},~~~ l=1,2,\cdots,n, \end{equation}
 is a set of solutions to the system \eqref{eq:EV}. Since
  $\textbf{A}^T$ is  positive definite, we deduce that the expression given in  \eqref{eq:solEV}  uniquely solves  \eqref{eq:EV}. Therefore, the optimal  reinsurance-investment  strategy   is given as follows
\begin{equation}\label{eq:opt_nH}\mH^*_l(u, x)=\frac{2\Lambda}{2\Lambda-\frac{\xi_{uu}}{\xi_u}}x,~~~~~l=1,2,\cdots,n,\end{equation} and \begin{equation} \label{eq:opt_npi}\pi^*(u)=-\displaystyle\frac{u-r}{\sigma^2}\frac{\xi_u}{\xi_{uu}}.\end{equation}
 Recall that $\Delta=\frac{(\mu-r)^2}{2\sigma^2}$. Bringing  \eqref{eq:opt_nH} and \eqref{eq:opt_npi} back into  the  differential operator $\mA^{\nu} \xi(u,m)$ and putting $\mA^{\nu} \xi(u,m)=0$, one can show that 
\begin{equation}\label{eq:ratio-nvar}ru-\kappa-\displaystyle \Delta\cdot\frac{\xi_u}{\xi_{uu}}+{A}_2\cdot\frac{2\Lambda^2}{2\Lambda-\frac{\xi_{uu}}{\xi_u}}=0\end{equation}
%\begin{equation}\label{eq:ratio4}\frac{\xi_{uu}}{\xi_u}=\frac{2\Lambda (ru-\kappa)+B+2\Lambda^2\bar{A}+\sqrt{(2\Lambda (ru-\kappa)+B+2\Lambda^2\bar{A})^2-8\Lambda B(ru-\kappa)}}{2(ru-\kappa)}\end{equation}
with 
\begin{equation}\label{eq:A2}
{A}_2=\sum_{l=1}^n(\la+\la_l)\E (X^{(l)})^2+\la\sum_{k\neq j}\E\big(X_i^{(k)}\big)\E\big(X_i^{(j)}\big).
\end{equation}
 Therefore, we have the following theorem.
 
 %Following the arguments and steps in Theorems \ref{th3.1} and  \ref{th4.1},  we can calculate the probability of drawdown  with the corresponding optimal reinsurance-investment strategy  in  following theorem.
 \begin{theorem}\label{thm:5.2}
 Let $\bar{g}_2(u,m)$and $\bar{f}_2(y)$ be given in Appendix  A.2. The minimum probability of drawdown for $n$ dependent classes of insurance business on $\mO$ is given by
$$ \phi(u,m)=\left\{\begin{array}{llll}
 \displaystyle1-\frac{\bar{g}_2(u,m)}{\bar{g}_2(u_s,m)},& \alpha m\leq u \leq u_s\leq m, \\[5mm]
 \displaystyle1- \exp\left\{-\int_m^{u_s}{\bar{f}_2(y)dy}\right\}\cdot\frac{\bar{g}_2(u,m)}{\bar{g}_2(u_s,u_s)}, & \alpha m\leq u \leq m \leq u_s.
 \end{array}\right.
$$
Also,  the corresponding  optimal reinsurance-investment strategy, for all $t\geq0$ and $l=1,2,\dots,n,$ is given in feedback form by
\begin{equation}\label{eq:strategy5}
 \mH^*_{lt}=\displaystyle\frac{2\Lambda}{2\Lambda-{\chi}_2(\hat{U}_t)}X^{(l)},~~~~~\pi^*_t=\displaystyle\frac{\mu-r}{\sigma^2}\cdot\frac{1}{{\chi}_2(\hat{U}_t)},
\end{equation}where the function $\chi_2(u)$ denotes the expression of $\chi_1(u)$ in \eqref{eq:beta1}, but replaces   $A_1$  by $A_2$.
\end{theorem}
\begin{remark}\label{rem:5.2}
Setting $l=n+1$, $p_{lj}=1$ for $j=1,2,\dots,n$, $p_{ij}=0$ $(i\neq j)$ for $i,j=1,2,\dots,n$, and $p_{ii}=1$ for $i=1,2,\dots,n,$ in Han et al. \cite{HLY18} (where the authors   investigated  optimal proportional reinsurance problem in a risk model with the thinning-dependence structure for both of the expected value principle and the variance premium principle), then the resulting  thinning-dependence risk model  reduces to the one  with common-shock.  We have proved that under the variance premium principle, the optimal retention strategy is in the form of quota-share reinsurance, which we can see clearly from \eqref{eq:strategy3} and \eqref{eq:strategy5}. Therefore,   setting $\pi=0$ in \eqref{eq:Uhat}, it is to be expected that the optimal results are consistent with   those in Han et al. \cite{HLY18}.
\end{remark}
\section{Some properties and numerical examples}\label{sec:6}
In this section, we first study some properties of the optimal  reinsurance-investment strategy and the corresponding value function in Corollary \ref{cor:6.1}. We use the expressions of $(\mH^*_1,\mH^*_2,\pi^*)$ and  $\phi$ given in  \eqref{eq:drawdown1}-\eqref{eq:strategy2}  for the  expected value principle and \eqref{eq:drawdown5}-\eqref{eq:beta1} for the variance premium principle, but in the interest of space, we omit the proof. When we write ``increases" or  `` decreases", we mean in the weak, or non-strict, sense. 

\begin{corollary}\label{cor:6.1} The optimal  reinsurance-investment strategy  and the corresponding value function satisfy the following properties:

\noindent
(i) The optimal strategies $\mH^*_i(u,x)$ $(i=1,2)$ and  $\pi^*(u)$   are decreasing functions with respect  to $u$, and tend to $0$ as the surplus $u$ approaches  $u_s$;

\noindent
(ii) $\mH^*_i(u,x)$ and  $x-\mH^*_i(u,x)$ $(i=1,2)$  both increase as the claim $x$ increases;

\noindent 
(iii) Under the variance premium principle, the optimal retention strategy $\mH^*_i$ $(i=1,2)$ increases with respect to $\Lambda$ and $\la$;

\noindent
(iv) The optimal investment strategy $\pi^*$ decreases as $\sigma$ increases;

\noindent
(v) The value function $\phi(u,m)$  decreases with respect to $u$, but decreases with respect to $m$ for both the cases of $m\geq u_s$ and $m<u_s$.

\end{corollary} 

Here,  in contrast to the variance premium principle,  the monotonicity of the optimal strategy  with respect to $\eta_i$ $(i=1,2)$ and $\la$ is not  clear  under the expected value  principle. This is due to the complex expression of  function $G$ given in \eqref{eq:G}.  Thus, in the following context, we   present  a numerical example to show  the monotonicity of the  function $G(7,d)$ with respect to $d$, and   then  carry on some numerical examples to investigate the effect of $u$, $\eta_1$  and $\la$  on the optimal reinsurance and investment  strategies.%Moreover, we find that the optimal retention and investment  strategies are decreasing functions with respect to the surplus value $u$, which is consistent  with the one in Corollary \ref{cor:6.1}.
   
   We assume that the insurer has two lines of business: one is heavy-tailed risk with small arriving intensity, and the other is light-tailed risk with large arriving intensity.  Note that the dependence between the two lines of business is determined by the intensity parameter $\la$. 
%The insurer and the reinsurer charge more premiums for the heavy-tailed risk business. 
 Let 
$$F_X(x)=1-\frac{1}{(x+1)^3},~~x\geq0;~~~F_Y(y)=1-e^{-2y},~~y\geq0.$$
Then we have $E(X)=\frac{1}{2}$, $E(Y)=\frac{1}{2}$. % $E(X^2)=1$, $E(Y^2)=\frac{1}{2}$. 
 Here, we assume that the insurer's premium is also calculated by  the expected value  principle, that is,
$$c=(1+\theta_1)(\la+\la_1)\E X+(1+\theta_2)(\la+\la_2)\E Y,$$
where $\theta_i(i = 1, 2)$ is  the insurer's  safety loading of the two classes of the insurance business, respectively.

\begin{example}\label{exm:6.1}We set $r=0.1$, $\eta_1=0.4$, $\eta_2=0.25$, $\theta_1=0.2$, $\theta_2=0.1$, $\la=2$, $\la_1=3$, $\la_2=5$, which implies that the safe level $u_s=10.25$. Also, we set $\alpha=0.1$ and $m=20$, which implies that the drawdown level equals 2. Finally, we set $\mu=0.5$ and $\sigma=1$.\end{example}
\begin{figure}[h]
\begin{tabular}{cc}
\begin{minipage}[t]{0.5\linewidth}
{\includegraphics[width=2.2 in]{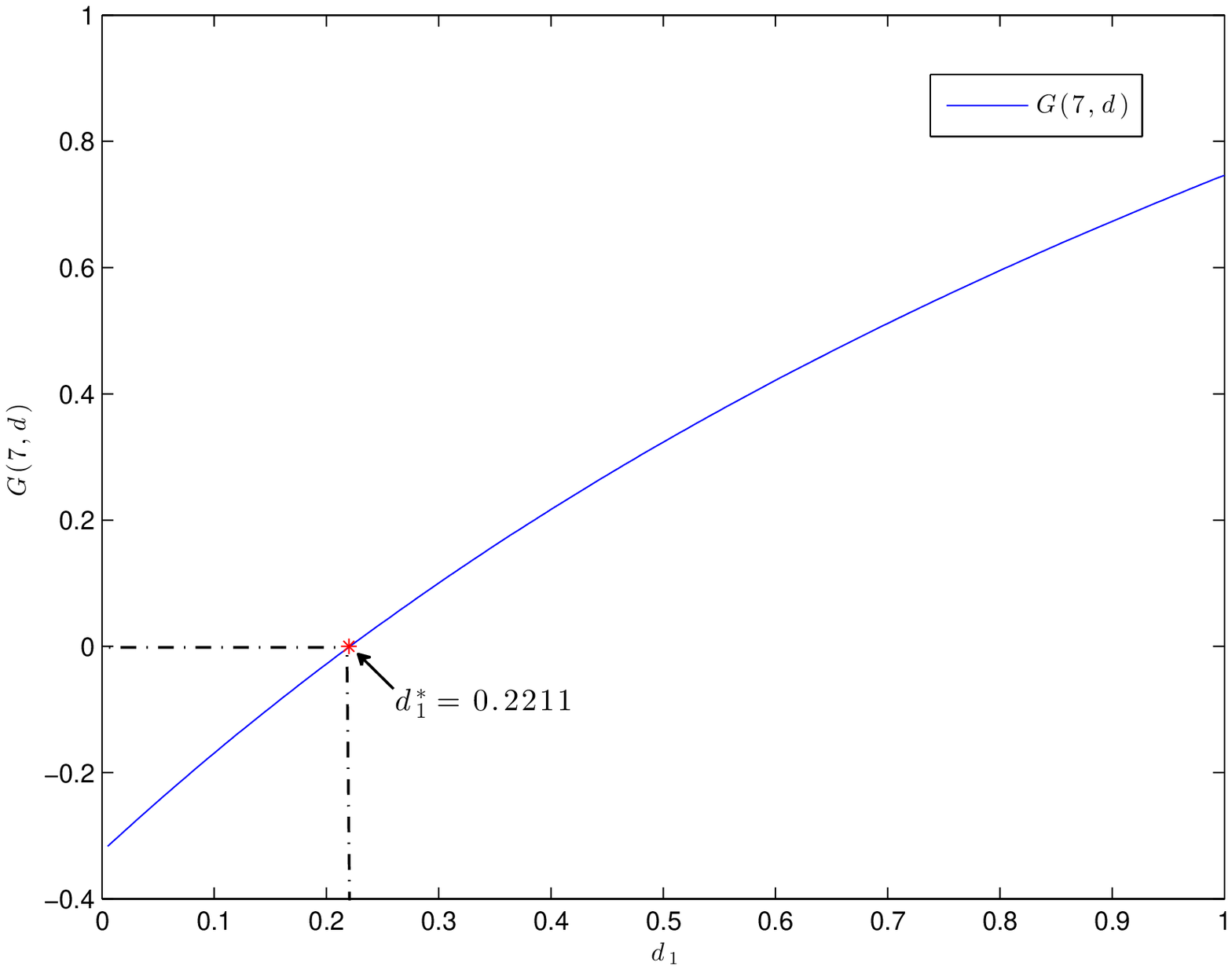}}
 \caption{The function  $G(7,d)$ }
\end{minipage}
\begin{minipage}[t]{0.5\linewidth}
{\includegraphics[width=2.2 in]{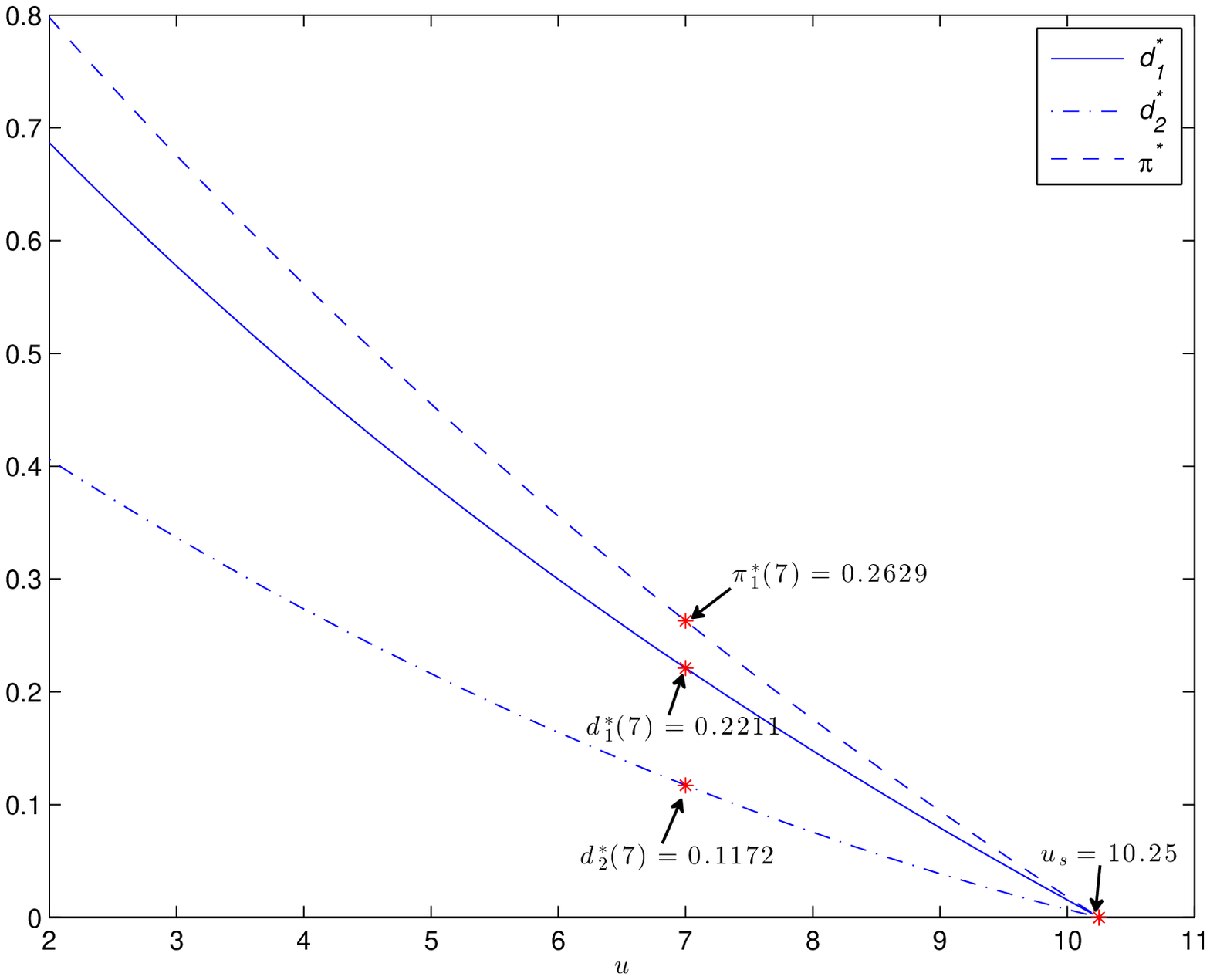}}
 \caption{The effect of $u$ on optimal strategy}
\end{minipage}\\
\end{tabular}
%\vspace{-1 in}
\end{figure}

We can see from Figure $1$ that the function $G(7,d)$ given in \eqref{eq:G} is an increasing function with respect to $d$.  Note that  we have $a_l=0$ with the given parameters. Thus, according to Lemma \ref{lem:4.1},  it is to be expected that there  exists a unique $d^*_1(u) \in[0,\infty)$ such that  $G(u, d^*_1)=0$.  For $u=7$, we obtain $d^*_1=0.2211$. Figure $2$ shows that the optimal  reinsurance  and investment  strategies  decrease as $u$ increases. These observations are kind of reasonable. When the value of the surplus increases toward $u_s$, the insurer can buy full reinsurance and invest all the surplus in the risk-free asset to earn interest rate r. As a result, the surplus will never drop below its current value, and drawdown cannot happen. In particular, we can see from Figure $2$ that  we have $d^*_1=0.2211$ for $u=7$, which is identical  to the one we find in Figure $1$.

\begin{example}\label{exm:6.2} We set $u=7$, $r=0.1$,  $\eta_2=0.25$, $\theta_1=0.2$, $\theta_2=0.1$, $\la_1=3$, $\la_2=5$. We set $\la=2$ and $\eta_1\in[0.3,0.5]$ for Figure $3$, and $\eta_1=0.4$ and $\la\in[1,6]$ for Figure $4$. Also, we set $\alpha=0.1$ and $m=20$, which implies that the drawdown level equals 2. Finally, we set $\mu=0.5$ and $\sigma=1$.\end{example}
\begin{figure}[h]
\begin{tabular}{cc}
\begin{minipage}[t]{0.5\linewidth}
{\includegraphics[width=2.2 in]{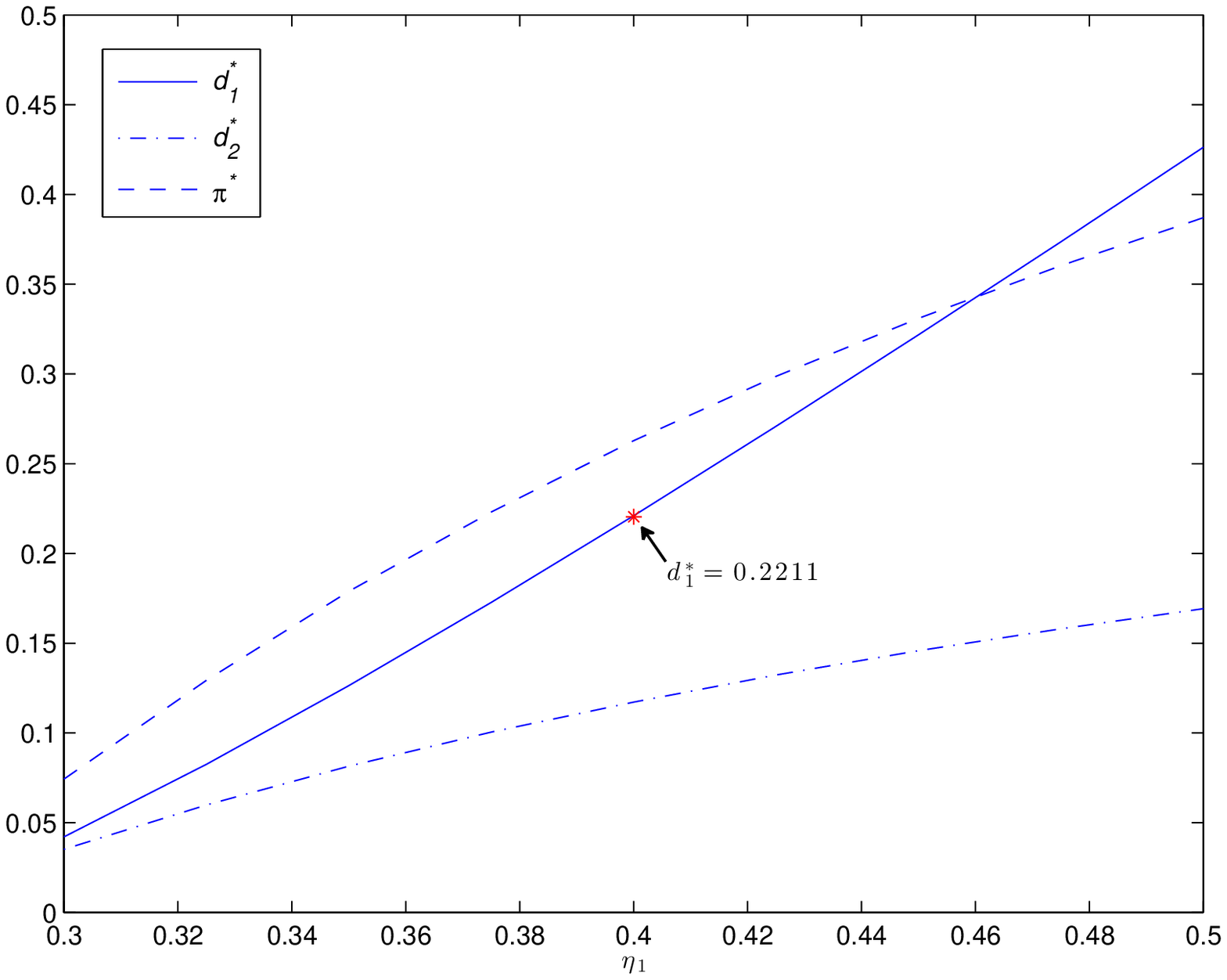}}
 \caption{The effect of $\eta_1$ on optimal strategy }
\end{minipage}
\begin{minipage}[t]{0.5\linewidth}
{\includegraphics[width=2.2 in]{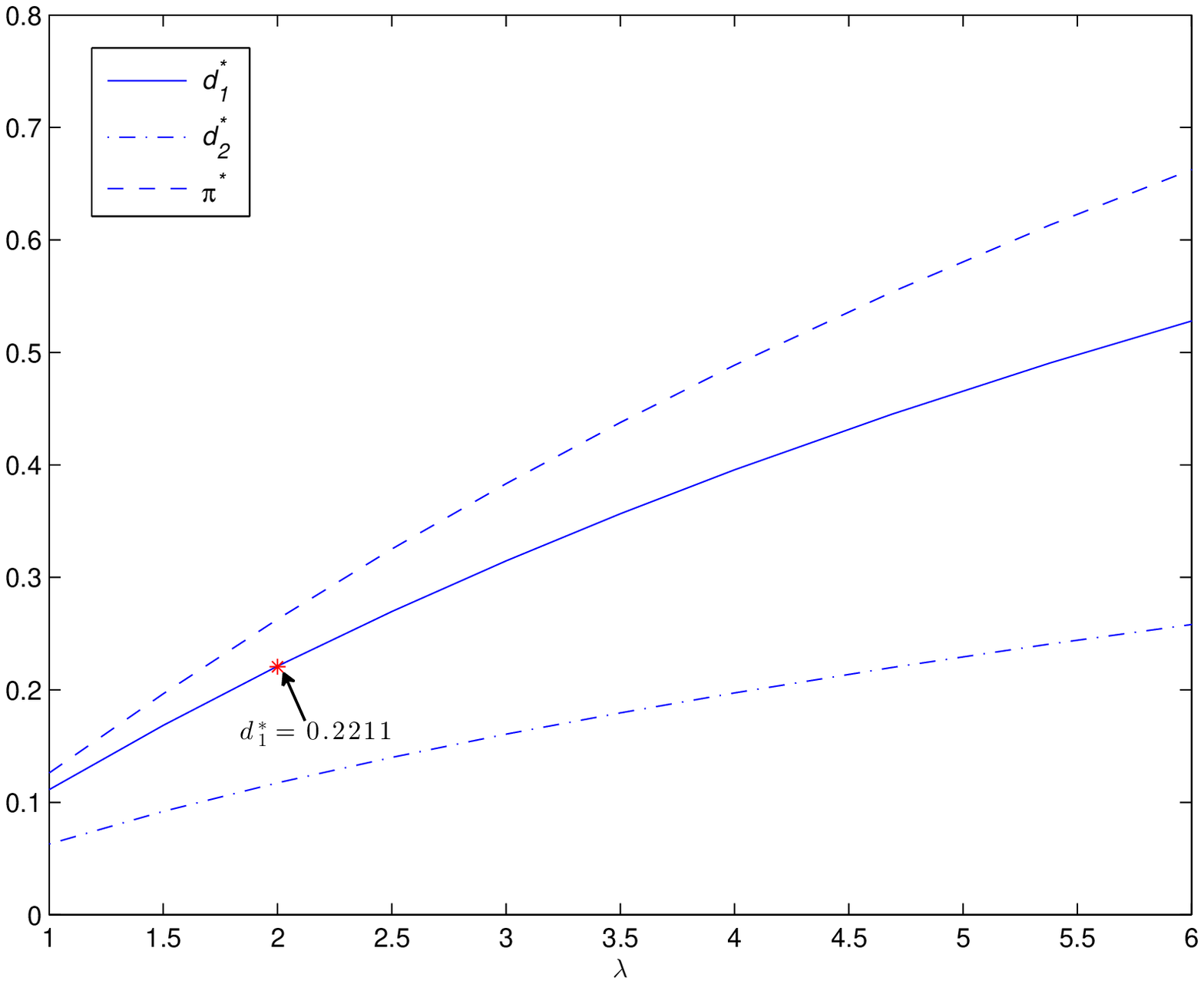}}
 \caption{The effect of $\la$ on optimal strategy}
\end{minipage}\\
\end{tabular}
%\vspace{-1 in}
\end{figure} 

In this example, we investigate the impact of the parameters of $\eta_1$ and $\la$ on the optimal strategy under the expected value principle. It shows that the reinsurance and investment strategies  increase as $\eta_1$ and $\la$ increase,  and the monotonicity is similar to the one under the variance premium principle.  This phenomenon  illustrates the intuition that when $\eta_1$ and $\la$ increase,  the reinsurance becomes more expensive, and then the insurer optimally purchases less reinsurance. Even though a larger value of $\la$  leads to more claims on average, %so one might think that the insurer would retain less of each claim, 
having to pay more in reinsurance premium dominates the insurer's optimal decision. Besides, we point out that the insurers always try to keep a greater share of each claim when the net premium is high $(\eta_1-\theta_1>\eta_2-\theta_2)$. When  the reinsurance premium keeps increasing, to avoid drawdown, the insurer might optimally invest larger amount in the risky asset to increase its profit. Note that the safe level also increases as $\eta_1$ and $\la$ increase. Specifically, the value of $d^*_1$ equals to the one given in Figure $1$ when we fix $\eta_1=0.4$ for Figure $3$ and $\la=2$ for Figure $4$, respectively.

\section{Conclusions}\label{sec:7}
In this paper, from an insurer's point of view, we find the optimal  reinsurance-investment strategy under a criterion of minimizing the probability of drawdown in a common shock dependence structure.   We assume that the surplus process of the insurer follows a diffusion, specifically, the diffusion approximation of the classical Cram\'er-Lundberg model.   The insurer can purchase per-loss  reinsurance and invest its surplus in a financial market consisting of one risky asset and one risk-free asset; thus, a finite safe level exists for our problem.  We first prove a verification theorem,  which is used to find the value function.  Then under two different premium principles, the explicit expressions for the optimal strategy and the associated value function are obtained which strongly depend on the  values of the wealth and the drawdown level. Specifically,  we observe  that the optimal reinsurance strategy is in the form of pure excess-of-loss reinsurance strategy under the expected value principle, and under the variance premium principle, the optimal reinsurance strategy is in the form of  pure quota-share reinsurance. In addition, we  extend our work to the  case where the  insurance company involves $n$ $(n\geq3)$ dependent classes of insurance business, and the optimal results are also derived explicitly.

Although the literature on optimal reinsurance  is increasing rapidly, there are still many interesting problems that deserve to be investigated.   For example,  we may study the optimal reinsurance-investment strategy which aims to minimize the probability that drawdown occurs over the time horizon $[0,e_\lambda]$ (where $e_{\lambda}$ is an independent exponentially distributed random time). Due to the presence of the time factor $e_\lambda$,   the corresponding HJB equation appears difficult to solve analytically.  Besides, most researchers only focus on the wealth management of an insurer and ignore the interest of a reinsurer. Actually, the reinsurer also wants to minimize the probability of drawdown. Thus, one may consider the optimization problem for a general insurance company which holds shares of an insurance company and a reinsurance company in a continuous model. Both of these problems are our future research directions.

\vspace{0.5cm}
\noindent{\large {\bf Acknowledgements}}
\vspace{0.3cm}

\noindent
This research was supported by the National Natural Science Foundation of China (Grant No. 11471165).

  \begin{appendices}
\section{Auxiliary functions}
\textbf{A.1.} The functions $g_{i}$ $(i=1,2,3)$ is  given by
$$\left\{\begin{array}{ll}
g_{1}(u,m)
=\displaystyle\int_{\alpha m}^u \exp\left\{-\int_{\alpha m}^y\frac{(\la+\la_1)\eta_1}{(\la+\la_1)d^*_1(w)+\la h_Y(d^*_2(w))}
dw\right\}dy, \\[5mm]
g_{2}(u,m)
=\displaystyle\int_{\alpha m}^{ \widetilde{u}}\exp\left\{-\int_{\alpha m}^y \frac{(\la+\la_1)\eta_1}{(\la+\la_1)d^*_1(w)+\la h_Y(d^*_2(w))}\right\}dy \\[5mm]
\quad\quad\quad\quad\quad+\displaystyle\int_{\widetilde{u}}^u \exp\left\{-\left(\int_{\alpha m}^{\widetilde{u}} \frac{(\la+\la_1)\eta_1}{(\la+\la_1)d^*_1(w)+\la h_Y(d^*_2(w))}+\int_{\widetilde{u}}^y\frac{\eta_1}{a_k(w)}\right)dw\right\}dy,\\[5mm]
g_3(u,m)
=\displaystyle\int_{\alpha m}^u \exp\left\{-\int_{\alpha m}^y\frac{\bar\Delta}{rw-\kappa}dw\right\}dy;
\end{array}\right.$$
and the functions $f_{i}$ $(i=1,2,3)$ is given by 
$$\left\{\begin{array}{ll}f_{1}(y)=\alpha \left[\displaystyle\frac{1}{g_{1}(y,y)}-\frac{(\la+\la_1)\eta_1}{(\la+\la_1)d^*_1(\alpha y)+\la h_Y(d^*_2(\alpha y))}\right],\\[5mm]
f_{2}(y)=\alpha \left[\displaystyle\frac{1}{g_{2}(y,y)}-\frac{(\la+\la_1)\eta_1}{(\la+\la_1)d^*_1(\alpha y)+\la h_Y(d^*_2(\alpha y))}\right],\\[5mm]
f_3(y)=\alpha \left[\displaystyle\frac{1}{g_3(y,y)}-\frac{\bar\Delta}{r\alpha y-\kappa}\right]\end{array}\right.$$
\textbf{A.2.} The functions $\bar{g}_i$ $(i=1,2)$ are given by
$$
\bar{g}_i(u,m)
=\displaystyle\int_{\alpha m}^u \exp\left\{-\int_{\alpha m}^y\chi_i(w)
dw\right\}dy;$$
and the functions $\bar{f}_i$ (i=1,2) is given by 
$$
 \bar{f}_i(y)=\alpha \left[\displaystyle\frac{1}{\bar{g}_i(y,y)}-\chi_i(\alpha y)\right].$$
\end{appendices}

\begin{thebibliography}{}
	%\small \setlength{\itemsep}{-.8mm}
\bibitem{ABY16a} Angoshtari, Bahman, Erhan Bayraktar, and Virginia R. Young (2016a). Optimal investment to minimize the probability of drawdown. {\it Stochastics}, \textbf{88}(6), 946-958.

\bibitem{ABY16b} Angoshtari, Bahman, Erhan Bayraktar, and Virginia R. Young (2016b).  Minimizing the probability of lifetime drawdown under constant consumption. {\it Insurance: Mathematics and Economics}, \textbf{69}, 210-223.

\bibitem{BCZ13} Bai, Lihua,  Jun Cai, and   Ming Zhou (2013). Optimal reinsurance policies for an insurer with a bivariate reserve risk process in a dynamic setting. {\it  Insurance: Mathematics and Economics}, \textbf{53}(3), 664-670.

\bibitem{BZ15} Bayraktar, Erhan and Yuchong Zhang (2015). Minimizing the probability of lifetime ruin under ambiguity aversion. {\it SIAM Journal on Control and Optimization}, \textbf{53}(1), 58-90.

\bibitem{B60} Borch, Karl (1960). An attempt to determine the optimum amount of stop loss reinsurance.{\it In: Transaction of the $16$th International Congress of Actuaries}, pp. 597-610.

\bibitem{CLLL15} Chen, Xinfu, David Landriault, Bin Li, and Dongchen Li (2015). On minimizing drawdown risks of lifetime investments. {\it Insurance: Mathematics and Economics}, \textbf{65}, 46-54.

\bibitem{CK95} Cvitani\'{c}, Jaksa and Ioannis Karatzas (1995). On portfolio optimization under drawdown constraints. {\it IMA Lecture Notes in Mathematical Applications}, \textbf{65}, 77-88.

%\bibitem{G91} Grandell, Jan (1991). Aspects of Risk Theory. Springer, New York.

\bibitem{GZ93} Grossman, Sanford J. and Zhongquan Zhou (1993). Optimal investment strategies for controlling drawdowns. {\it Mathematical Finance}, \textbf{3}(3), 241-276.

\bibitem{HLY19} Han,  Xia, Zhibin Liang, and Virginia R. Young (2020). Optimal reinsurance strategy to minimize the probability of drawdown under a Mean-Variance premium principle.  {\it Scandinavian Actuarial Journal}. DOI:10.1080/03461238.2020.1788136..

\bibitem{HLY18} Han, Xia,  Zhibin Liang, and  Kam Chuen Yuen (2018). Optimal proportional reinsurance to minimize the probability of drawdown under thinning-dependence structure. {\it Scandinavian Actuarial Journal},  \textbf{2018}(10), 863-889.



\bibitem {HT10} Hipp, Christian and Michael  Taksar (2010). Optimal non-proportional reinsurance control. {\it Insurance: Mathematics and  Economics}, \textbf{47}(2), 246-254.


\bibitem{KOP17} Kardaras, Constantinos, Jan Ob\l\'oj, and Eckhard Platen (2017). The Num\'eraire Property and Long-Term Growth Optimality for Drawdown-Constrained Investments. { \it Mathematical Finance}, \textbf {27}(1), 68-95.

\bibitem{LRZ17}  Li, Danping, Ximin Rong, and Hui Zhao (2017). Equilibrium excess-of-loss reinsurance-investment strategy for a mean-variance insurer under stochastic volatility model. {\it Communications in Statistics-Theory and Methods}, \textbf{46}(19), 9459-9475.

\bibitem{LY18} Liang, Xiaoqing and Virginia R.Young (2018). Minimizing the probability of ruin: optimal per-loss reinsurance.   {\it Insurance: Mathematics and  Economics}, \textbf{82}, 181-190.

 \bibitem{LG11} Liang, Zhibin and Junyi Guo (2011). Optimal combining quota-share and excess of loss reinsurance to maximize the expected utility. {\it Journal of Applied Mathematics and Computing}, \textbf{36} (1-2), 11-25.
 


 \bibitem{LY16} Liang, Zhibin and Kam Chuen Yuen (2016). Optimal dynamic reinsurance with dependence risks: variance premium principle. {\it Scandinavian Actuarial Journal},  \textbf{2016}(1), 18-36.
 
%\bibitem{LWZ19} Luo, Shangzhen, Mingming Wang, and Wei Zhu (2019). Maximizing a robust goal-reaching probability with penalization on ambiguity. Journal of Computational and Applied Mathematics, \textbf{348}, 261-281.

 \bibitem{PS85} Pestien, Victor C. and William D. Sudderth (1985). Continuous-time red and black: how to control a diffusion to a goal. {\it Mathematics of Operations Research}, \textbf{8}(2),  11-31.

\bibitem{PY05} Promislow, S. David  and Virginia R. Young (2005). Minimizing the probability of ruin when claims follow Brownian motion with drift. {\it North American Actuarial Journal}, \textbf{9}(3), 110-128.

%\bibitem{S02} Schmidli, Hanspeter (2002). On minimizing the ruin probability by investment and reinsurance. {\it The Annals of Applied Probability}, \textbf{12} (3), 890-907.

\bibitem{WY05} Wang, Guojing and  Kam Chuen Yuen (2005). On a correlated aggregate claims model with thinning-dependence structure. {\it  Insurance: Mathematics and Economics}, \textbf{36}, 456-468.


\bibitem{YGW02}  Yuen, Kam Chuen, Junyi Guo,  and Xueyuan Wu (2002). On a correlated aggregate claims model with Poisson and Erlang risk process. {\it Insurance: Mathematics and Economics}, \textbf{31}, 205-214.
%\bibitem{WRZ 18} Wang, Yajie, Ximin Rong, and Hui Zhao (2018).  Optimal investment strategies for an insurer and a reinsurer with a jump diffusion risk process under the CEV model. {\it Journal of Computational and Applied Mathematics},  \textbf{328},  414-431.
\bibitem{ZMZ16} Zhang, Xin, Hui Meng, and Yan Zeng (2016). Optimal investment and reinsurance strategies for insurers with generalized mean-variance premium principle and no-short selling. {\it  Insurance: Mathematics and Economics}, \textbf{67}, 125-132.

	\end{thebibliography}


\end{document}